\documentclass[a4paper,11pt]{article}
\usepackage{ifstyle}   
\usepackage{graphicx}
\usepackage{epsfig} 
\usepackage{color}
\numberwithin{equation}{section}  

\usepackage[francais]{babel}               
\usepackage[T1]{fontenc}
\usepackage[latin1]{inputenc}  
\usepackage{amssymb,diagrams}             
\usepackage{array}                 
\usepackage{fullpage}

\newtheorem{lem}[subsubsection]{Lemme}
\newtheorem{iw}[subsection]{Lemme de Wirtinger} 
\newtheorem{mut}{Conjecture}

\newcommand{\N}{\mathbb N}
\newcommand{\Z}{\mathbb Z}  
\newcommand{\Q}{\mathbb Q}
\newcommand{\R}{\mathbb R}
\newcommand{\C}{\mathbb C}

\newcommand{\I}{\mathbb I}
\newcommand{\proj}{\mathbb P}
\newcommand{\contract}{\mathrel{\kern-1.5pt\vrule width6.0pt height0.4pt depth0pt
                \vrule width0.4pt height4.0pt depth0pt}}
\newcommand{\retract}{\mathrel{\kern-1.5pt \vrule width0.4pt height5.0pt depth0pt
                \vrule width6.0pt height0.4pt depth0pt }}

\def\opn#1#2{\def#1{\operatorname{#2} } } 
\opn\Ric{Ricci} 
\opn\Trac{Trace} 
\opn\det{det} 
\opn\Ker{Ker} 
\opn\exp{exp}
\opn\exph{exph}
\opn\Herm{Herm} 
 
\begin{document} 
\begin{center} 
\huge{\bf{Fonctions plurisousharmoniques et courants positifs de type $(1,1)$ sur une variété presque complexe}}
\\
\vspace{0.4cm}
\Large{Nefton Pali}
\\
\vspace{0.4cm}
\large{Institut Fourier, UMR 5582, Université Joseph Fourier
\\
BP 74, 38402 St-Martin-d'Hères cedex, France}
\\
\large{E-mail: \textit{nefton.pali@ujf-grenoble.fr}}
\end{center} 
{\bf Résumé}.-Une fonction semi-continue supérieurement $u$ sur une variété presque complexe $(X,J)$ est dite plurisousharmonique si la restriction à toute courbe pseudo-holomorphe locale est sous-harmonique. 
Comme dans le cas analytique complexe, nous conjecturons que la notion de plurisousharmonicité pour une fonction $u$ est équivalente à la positivité du $(1,1)$-courant $i\partial_{_J }\bar{\partial}_{_J }u$, (lequel n'est pas forcément fermé dans le cas non intégrable). La conjecture est triviale dans le cas d'une fonction $u$ de classe ${\cal C}^2$. 
Le résultat en question est élémentaire dans le cas complexe intégrable car l'opérateur $i\partial_{_J }\bar{\partial}_{_J }$ s'écrit comme un opérateur à coefficients constants dans des coordonnées complexes. On peut donc facilement conserver la positivité du courant en régularisant avec des noyaux $\ci$ usuels. Dans le cas presque complexe non intégrable ceci ce n'est pas possible et la preuve du résultat exige un étude beaucoup plus intrinsèque. 
Nous montrons la nécessité de la positivité du $(1,1)$-courant $i\partial_{_J }\bar{\partial}_{_J }u$ en utilisant la théorie locale des courbes $J$-holomorphes. Nous montrons aussi  la suffisance de la positivité dans le cas particulier d'une fonction $f$ semi-continue supérieurement et continue en dehors du lieu singulier $f^{-1}(-\infty)$.
Pour prouver la suffisance de la positivité dans le cas général où $u$ est une distribution réelle nous proposons une méthode qui utilise un délicat argument de régularisation des fonctions plurisousharmoniques introduit par Demailly (\cite{Dem-2}). La méthode consiste à régulariser la fonction $u$ à l'aide du flot géodésique induit par une connexion de Chern sur le fibré tangent de la variété presque complexe, (cf. \cite{Pal}).
\\
\\
{\bf{Abstract}}.-If $(X,J)$ is an almost complex manifold, then a function $u$ is said to
be plurisubharmonic on $X$ if it is upper semi-continuous and its
restriction to every local pseudo-holomorphic curve is subharmonic. As
in the complex case, it is conjectured that plurisubharmonicity is
equivalent to the fact that the 
$(1,1)$-current $i\partial_{_J }\bar{\partial}_{_J }u$ is positive, (the 
$(1,1)$-current $i\partial_{_J }\bar{\partial}_{_J }u$ need not be closed here!). The conjecture is trivial if $u$ is of class ${\cal C}^2$. The result is elementary in the complex integrable case because the operator $i\partial_{_J }\bar{\partial}_{_J }$ can be written as an operator with constant coefficients in complex coordinates. Hence the positivity of the current is preserved by regularising with usual convolution kernels. This is not possible in the almost complex non integrable case and the proof of the result requires a much more intrinsic study. In this chapter we prove the necessity of the positivity of the $(1,1)$-current $i\partial_{_J }\bar{\partial}_{_J }u$. We prove also the sufficiency of the positivity in the particular case of an upper semi-continuous function $f$ which is continuous in the complement of the singular locus $f^{-1}(-\infty)$.
For the proof of the sufficiency of the positivity in the general case of a real distribution $u$, we suggest a method depending on a rather delicate regularisation argument introduced by Demailly (\cite{Dem-2}). This method consists of regularing the function $u$ by means of the flow induced by a Chern connection on the tangent bundle of the almost complex manifold, (see \cite{Pal}).
\tableofcontents
\section{Préliminaires}
Soit $(X,J)$ une variété presque complexe de classe ${\cal C}^{\infty}$ et de dimension réelle $2n$. On désigne par ${\cal E}_X\equiv{\cal E}_X(\R)$ le faisceau des fonctions ${\cal C}^{\infty}$ à valeurs réelles, par $\pi^{1,0} _{_J}:T_X\otimes_{_{\R}}\C\longrightarrow T^{1,0}_{X,J}$ la projection sur le fibré des $(1,0)$-vecteurs tangents et par $\pi^{0,1} _{_J}$ celle sur le fibré des $(0,1)$-vecteurs tangents. On désigne par $T_{X,J}$ le fibré tangent dont les fibres sont munies de la structure complexe donnée par $J$ et par
$$
{\cal E}^{p,q}_{_{X,J} }\equiv {\cal E}(\Lambda ^{p,q}_{_J}T_X^*),
\;\Lambda ^{p,q}_{_J}T_X^*:=\Lambda^p _{_{\C}}(T^{1,0}_{X,J}  )^*\otimes_{_{\C}} \Lambda^q _{_{\C}}(T^{0,1}_{X,J})^*
$$ 
le faisceau des $(p,q)$-formes par rapport à la structure presque complexe $J$. On rappelle que sur une variété presque complexe la différentielle se décompose sous la forme 
$$
d=\partial_{_J }+\bar{\partial}_{_J }-\theta_{_J }-\bar{\theta }_{_J },
$$ 
où pour toute $k$-forme complexe 
$\omega \in {\cal E}(\Lambda^k _{_{\C}}(T_X\otimes_{_{\R}}\C)^*)(U)$ au dessus d'un ouvert $U$ et tout champ de vecteurs complexes $\xi _0,...,\xi _k\in{\cal E}(T_X\otimes_{_{\R}}\C)(U)$ on a les expressions suivantes:
\begin{eqnarray*} 
&\displaystyle{\partial_{_J }\omega\, (\xi _0,...,\xi _k):=\sum_{0\leq j \leq k}(-1)^j\xi^{1,0}_j .\,\omega (\xi _0,...,\widehat{\xi _j},..., \xi _k)+ }&
\\
\\
&\displaystyle{+\sum_{0\leq j<l \leq k}(-1)^{j+l}\omega ([\xi^{1,0}_j,\xi^{1,0}_l]^{1,0}+[\xi^{0,1} _j,\xi^{1,0} _l]^{0,1}+[\xi^{1,0} _j,\xi^{0,1} _l]^{0,1},\xi _0,...,\widehat{\xi _j},...,\widehat{\xi _l},..., \xi _k)   }&
\\
\\
\\
&\displaystyle{\bar{\partial}_{_J }\omega\, (\xi _0,...,\xi _k):=\sum_{0\leq j \leq k}(-1)^j\xi^{0,1}_j .\,\omega (\xi _0,...,\widehat{\xi _j},..., \xi _k)+ }&
\\
\\
&\displaystyle{+\sum_{0\leq j<l \leq k}(-1)^{j+l}\omega ([\xi^{0,1}_j,\xi^{0,1}_l]^{0,1}+[\xi^{0,1} _j,\xi^{1,0} _l]^{1,0}+[\xi^{1,0} _j,\xi^{0,1} _l]^{1,0},\xi _0,...,\widehat{\xi _j},...,\widehat{\xi _l},..., \xi _k)   }&
\end{eqnarray*}
\begin{eqnarray*} 
&\displaystyle{\theta_{_J }\omega\, (\xi _0,...,\xi _k):=-\sum_{0\leq j<l \leq k}(-1)^{j+l}\omega ([\xi^{1,0}_j,\xi^{1,0}_l]^{0,1},\xi _0,...,\widehat{\xi _j},...,\widehat{\xi _l},..., \xi _k)   }&
\\
\\
\\
&\displaystyle{\bar{\theta }_{_J }\omega\, (\xi _0,...,\xi _k):=-\sum_{0\leq j<l \leq k}(-1)^{j+l}\omega ([\xi^{0,1}_j,\xi^{0,1}_l]^{1,0},\xi _0,...,\widehat{\xi _j},...,\widehat{\xi _l},..., \xi _k)   }&
\end{eqnarray*}
avec $\xi ^{1,0}:=\pi ^{1,0} _{_J}(\xi ),\;[\cdot,\cdot]^{1,0}:=\pi^{1,0}_{_J}[\cdot,\cdot] $ et de façon analogue pour les indices $(0,1)$. Les bidegrés des opérateurs 
$\partial_{_J },\,\bar{\partial}_{_J },\,\theta_{_J }$ et $\bar{\theta }_{_J }$ sont respectivement $(1,0),\;(0,1),\;(2,-1)$ et $(-1,2)$. En effet si 
$\omega\in{\cal E}^{p,q}_{_{X,J} } (U)$ est une $(p,q)$-forme alors les $(p+q+1)$-formes $\partial_{_J }\omega,\,\bar{\partial}_{_J }\omega,\,\theta_{_J }\omega,\,\bar{\theta }_{_J }\omega$ sont nulles en restriction aux fibrés $\Lambda ^{r,s}_{_J}T_X,\,r+s=p+q+1$ respectivement aux bidegrés $(r,s)\not=(p+1,q),\,(r,s)\not=(p,q+1),\,(r,s)\not=(p+2,q-1)
,\,(r,s)\not=(p-1,p+2)$. On déduit alors que l'opérateur $T=\partial_{_J },\;\bar{\partial}_{_J },\;\theta_{_J }$ où $\bar{\theta }_{_J }$ vérifie la règle de Leibnitz 
$$
T(u\wedge v)= Tu\wedge v+(-1)^{\deg\,u} u\wedge  Tv.
$$
On a aussi les formules $\overline{(\partial_{_J }\,u)}=\bar{\partial}_{_J } \,\bar{u},\,\overline{(\theta_{_J }\,u)}=\bar{\theta }_{_J }\,\bar{u}   $.
\begin{defi}
On désigne par $\tau_{_J}\in {\cal E}(\Lambda ^{2,0}_{_J}T_X^*\otimes_{_{\C}} T^{0,1}_{X,J})(X)$ le tenseur de la torsion de la structure presque complexe définie par la formule $\tau_{_J}(\xi ,\eta):= [\xi^{1,0},\eta^{1,0}]^{0,1}$ pour tout $\xi,\eta\in{\cal E}(T_X\otimes_{_{\R}}\C)(U)$, où $U\subset X$ désigne un ouvert quelconque. Le tenseur de  la structure presque complexe est dit intégrable si $\tau_{_J}=0$.
\end{defi}
On remarque que $\tau_{_J}=0$ si et seulement si $\theta_{_J}=0$, si et seulement si $d=\partial_{_J }+\bar{\partial}_{_J }$. On rappelle le célèbre théorème de Newlander-Nirenberg (voir \cite{We}, \cite{Hor}, \cite{Dem-1}, chapitre VIII, \cite{mal}, \cite{Nij-Woo} et \cite{New-Nir}).
\begin{theo}{\bf (Newlander-Nirenberg)}.
Soit $(X,J)$ une variété presque complexe. L'existence d'une structure holomorphe ${\cal O}_X$ sur la variété $X$ telle que la structure presque complexe associée $J_{{\cal O}_X } $ soit égale à $J$ est équivalente à l'intégrabilité de la structure presque complexe $J$.
\end{theo} 
On désigne par ${\cal H}^p_m$ la mesure de Hausdorff $p$-dimensionnelle dans $\R^m$ et par $\lambda$ la mesure de Lebesgue sur $\R^m$.  On désigne par $B_{r }(x)$ la boule ouverte  de $\R^m$ de centre l'origine et de rayon $r >0$ et par $S_r(x)$ la sphère de dimension $m-1$ dans $\R^m$ de centre l'origine et de rayon $r >0$. Soit $f$ une fonction Borel-mesurable et localement bornée sur un ouvert $U\subset\R^m$. Pour tout $\overline{B_{r }(x)}\subset U$ on définit les quantités
$$
\mu_{B}(f,x,r):=\frac{1}{\lambda(B_r(x))}   \int\limits_{B_r(x) }f\,d\lambda \quad \mbox{et} \quad
\mu _{S}(f,x,r):=\frac{1}{{\cal H}^{m-1}_{m}(S_r(x))}   \int\limits_{S_r(x) }f\,d{\cal H}^{m-1}_{m}. 
$$
On a la définition suivante (cf. \cite{Dem-1} pour plus de détails).
\begin{defi}
Une fonction $f:U\longrightarrow [-\infty,+\infty)$ semi-continue supérieurement est dite sous-harmonique si elle vérifie une des deux propriétés équivalentes suivantes:
\\
a) $f(x)\leq \mu_{B}(f,x,r)$ pour tout $\overline{B_{r }(x)}\subset U$;
\\
b) $f(x)\leq \mu_{S}(f,x,r)$ pour tout $S_{r }(x)\subset U$.
\end{defi}
Si $f\in{\cal C}^2(U,\R)$ alors on déduit d'après la deuxième identité de Green que $f$ est sous-harmonique si et seulement si $\Delta f\geq 0$.
 De façon générale on a le théorème classique suivant (cf. \cite{Dem-1}, chapitre I).
\begin{theo}\label{equiv-SH}
Soit $f$ une fonction sous-harmonique sur un ouvert connexe $U$. Alors soit $f\equiv -\infty$, soit $f\in L^1_{loc}(U)$ et dans ce cas le Laplacien au sens des distributions $\Delta f$ est une mesure positive. Réciproquement soit  $u$ une distribution sur $U$ telle que $\Delta u$ soit une mesure positive. Alors il existe une unique fonction $f$ sous-harmonique sur $U$ telle que $u$ soit la distribution associée à $f$.
\end{theo}
On déduit d'après le théorème \ref{equiv-SH} q'une fonction $f$ est sous-harmonique sur un ouvert $U$ si et seulement si pour tout $x\in U$ il existe un voisinage ouvert $V_x\subset U$ de $x$ tel que $f$ est sous-harmonique sur $V_x$.
\\
On désigne par $j$ la structure presque complexe canonique sur $\R^2\equiv \C$, par $J_0$ la structure presque complexe canonique sur $\R^{2n} \equiv \C^n$ et par $B^1_{\delta} \subset\R^2$ la boule complexe de centre l'origine et de rayon $\delta$. On rappelle la définition suivante.
\begin{defi} 
Soit $(X,J)$ une variété presque complexe. Une courbe $J$-holomorphe locale est une application différentiable 
$\gamma :B^1_{\delta} \longrightarrow X$  telle que sa différentielle vérifie la condition 
$J(\gamma (z))\cdot d_z\gamma=d_z\gamma \cdot j$ pour tout $z\in B^1_{\delta }$.
\end{defi}  
On a la définition suivante.
\begin{defi} 
Soit $(X,J)$ une variété presque complexe. Une fonction $f:X\longrightarrow [-\infty,+\infty)$ semi-continue supérieurement est dite $J$-plurisousharmonique si pour toute courbe $J$-holomorphe locale $\gamma$ définie sur le disque $B_{\delta}\subset\R^2$, la 
composée $f\circ \gamma$ est sous-harmonique sur le disque $B_{\delta}$.
\end{defi} 
Nous désignerons par $Psh(X,J)$ l'ensemble des fonctions $J$-plurisousharmoniques. Si 
\\
$u\in {\cal D}'_{2n}(X,\C)$ est une distribution à valeurs complexes sur $X$ nous sommes particulièrement intéressés par le courant 
$i\partial_{_J }\bar{\partial}_{_J }u\in{{\cal D}'}^{1,1} (X)$. En général on désigne par ${{\cal D}'}^{k,k} (X) $ les sections globales du faisceau 
$$
{\cal E}(\Lambda ^{k,k}_{_J}\,T_X^*) \otimes_{_{{\cal E}_X(\C) }}{\cal D}'_{2n}(\C)
$$ 
où ${\cal D}'_{2n}(\C)$ représente le faisceau des distributions à valeurs complexes sur $X$. Il est bien connu que  ${{\cal D}'}^{k,k} (X)$ s'identifie naturellement par intégration au dual topologique ${\cal D}'_{n-k,n-k} (X)$ de l'espace ${\cal D}^{n-k,n-k} (X)$ des $(n-k,n-k)$-formes ${\cal C}^{\infty}$ à support compact muni de la topologie de la convergence localement uniforme de toutes les dérivées (cf. \cite{Dem-1}, chapitre I et \cite{DeR}). On utilisera pourtant dans la suite l'identification des notations 
${{\cal D}'}^{k,k} (X)\equiv {\cal D}'_{n-k,n-k} (X)$. Le courant $i\partial_{_J }\bar{\partial}_{_J }u$ s'écrit  explicitement sous la forme
\begin{eqnarray}\label{explicit(1,1)cour} 
i\partial_{_J }\bar{\partial}_{_J }u\,(\xi _0,\xi _1)=i\xi ^{1,0}_0 .\, \xi ^{0,1}_1.\,u-i\xi ^{1,0}_1 .\, \xi ^{0,1}_0.\,u-
i[\xi^{0,1} _0,\xi^{1,0} _1]^{0,1}.\,u-i[\xi^{1,0}_0,\xi^{0,1}_1]^{0,1}.\,u
\end{eqnarray}
pour tout champ de vecteurs complexes $\xi _0,\,\xi _1\in {\cal E}(T_X\otimes_{_{\R}}\C)(X)$. On rappelle que la dérivée $\xi .\,u$ d'une distribution $u$ par rapport à un champ de vecteurs $\xi$ est donnée par la formule 
$$
\langle\xi .\,u,\varphi \rangle :=-\langle u,d(\xi \contract\varphi )\rangle 
$$
pour tout  
$\varphi \in {\cal D}^{2n}(X,\C)$. On remarque que si la distribution $u$ est réelle alors le courant $i\partial_{_J }\bar{\partial}_{_J }u$ l'est aussi. En effet en degré zéro on a l'identité 
$\partial_{_J }\bar{\partial}_{_J }=-\bar{\partial}_{_J }\partial_{_J }$. Ceci découle de la relation 
$$
\partial_{_J }\bar{\partial}_{_J }+\bar{\partial}_{_J }\partial_{_J }=- \theta_{_J }  \bar{\theta }_{_J }-\bar{\theta }_{_J } \theta_{_J }
$$
et du fait que les opérateurs $\theta_{_J } $  et $\bar{\theta }_{_J } $ sont nuls en degré zéro.
Si $(\zeta _1,...,\zeta _n)$ est un repère local du fibré $ T^{1,0}_{X,J}$ alors l'expression du courant en question par rapport au repère choisi est:
\begin{eqnarray}\label{exploccour1,1}
i\partial_{_J }\bar{\partial}_{_J }u=i\sum_{1\leq k,l\leq n}(\zeta _k\,.\bar{\zeta}_l .\,u-[\zeta _k,\bar{\zeta}_l]^{0,1}.\,u )\,\zeta _k^*\wedge\bar{\zeta}_l^*
\end{eqnarray}
On remarque que dans le cas intégrable, si on considère un repère local  holomorphe $\zeta _k\in {\cal O}(T^{1,0}_{X,J})(U)$, $k=1,...,n$, on a $[\zeta _k,\bar{\zeta}_l]=0$ pour tout indice $k,l$. Rappelons maintenant la définition suivante:
\begin{defi} 
Une $(p,p)$-forme $u\in\Lambda ^{p,p}_{_J}T_{X,x}^*$ est dite positive si $u(\xi _1,J\xi _1,...,\xi _p,J\xi _p)\geq 0$ pour tout vecteur  
$\xi _1,...,\xi _p\in T_{X,x}$.
\\
Une $(q,q)$-forme  $v\in\Lambda ^{q,q}_{_J}T_{X,x}^*$ est dite fortement positive si elle peut être exprimée sous la forme 
$$
v=\sum \lambda _t\,i\alpha_{t,1}\wedge\bar{\alpha }_{t,1} \wedge...\wedge i\alpha_{t,q}\wedge\bar{\alpha }_{t,q}
$$
avec $\lambda _t\geq 0$ et 
$\alpha _{t,k}\in (T^{1,0}_{X,J,x}  )^*$
\end{defi} 
Bien évidemment l'ensemble des $(q,q)$-formes fortement positives est un cône convexe fermé. Il est bien connu que l'ensemble des $(p,p)$-formes  positives est le cône dual des $(q,q)$-formes fortement positives, où $q=n-p$, via la dualité donnée par le produit extérieur (cf. \cite{Dem-1}, chapitre III et \cite{Lel}). La dualité en question implique alors que toutes les formes positives sont réelles, (les formes fortement positives étant réelles). Soit 
$$
u=i^{p^2} \sum_{|K|=|H|=p}\,u_{K,H}\,\zeta _K^*\wedge\bar{\zeta}^*_H 
$$
une $(p,p)$-forme et $\xi_t =\sum_k \,(\lambda _{t,k} \,\zeta _k+\bar{\lambda}_{t,k}  \,\bar{\zeta}_k),\;t=1,...,p$,  vecteurs réels. On désigne par 
$\lambda =(\lambda _{t,k})\in M_{p,n}(\C)$ la $(p,n)$-matrice associée aux coefficients $\lambda _{t,k}$. Les identités
$$
\xi _1\wedge J\xi _1\wedge ...\wedge\xi _p\wedge J\xi _p=2^p(-i)^p\xi ^{1,0}_1\wedge \xi ^{0,1}_1\wedge ...\wedge \xi ^{1,0}_p\wedge \xi ^{0,1}_p
$$
et $i^{p^2}(-1)^{p(p-1)/2}=i^p$ impliquent les égalités
\begin{eqnarray*} 
u(\xi _1,J\xi _1,...,\xi _p,J\xi _p)=2^p(-i)^pu(\xi ^{1,0}_1,\xi ^{0,1}_1,...,\xi ^{1,0}_p,\xi ^{0,1}_p)=
\\
\\
=2^p(-i)^p(-1)^{p(p-1)/2} u(\xi ^{1,0}_1,...,\xi ^{1,0}_p,\xi ^{0,1}_1,...,\xi ^{0,1}_p)=\quad
\\
\\
=2^p \sum_{|K|=|H|=p}\,u_{K,H}\,\det\lambda _K\cdot \overline{(\det\lambda _H)} .\qquad\qquad
\end{eqnarray*}
On aura alors que la $(p,p)$-forme $u$ est positive si et seulement le dernier terme de l'égalité précédente est positif pour toute matrice $\lambda$. Dans le cas 
$p=1$ la matrice hermitienne $(u_{k,h})$ associée au coefficients de la forme $u$ est semidéfinie positive et une diagonalisation de celle ci montre qu'on peut exprimer $u$ sous la forme 
$
u=\sum_{1\leq t\leq r} i\,\alpha_t\wedge\bar{\alpha }_t ,
$
où $r$ est le rang de la forme $u$. On a donc que la notion de positivité coïncide avec celle de forte positivité en degré $(1,1)$ et par dualité aussi en degré $(n-1,n-1)$ (et bien évidement en bidegré $(0,0)$ et $(n,n)$). Nous montrons maintenant un premier résultat qui exprime la relation forte qui existe entre les formes positives et les fonctions plurisousharmoniques.
\begin{lem}
Soit $f\in {\cal C}^2(X,\R)$. Alors $f\in Psh(X,J)$ si et seulement si la forme $i\partial_{_J }\bar{\partial}_{_J}\,f$ est positive.
\end{lem} 
$Preuve$. Nous commençons par montrer la nécessité de la positivité de la forme $i\partial_{_J }\bar{\partial}_{_J}\,f$. En effet soit $\xi \in T_{X,x}$ un vecteur réel. Il existe alors une courbe $J$-holomorphe $\gamma$ telle que $\gamma(0)=x$ et $\xi =d\gamma\Big(\frac{\partial}{\partial x}{\vphantom{x} } _{| _{_0} } \Big)$ 
(voir par exemple l'article de Sikorav, théorème 3.1.1 dans l'ouvrage de Audin et Lafontaine \cite{Ad} ou la preuve du théorème $\ref{Jplat} $ qui suivra). Le fait que la courbe $\gamma$ soit $J$-holomorphe implique la première et troisième des égalités suivantes:
\begin{eqnarray*} 
i\partial_{_J }\bar{\partial}_{_J }f\,(\xi ,J\xi)=
i\partial_{_J }\bar{\partial}_{_J }f
\,\Big(d\gamma\Big(\frac{\partial}{\partial x}{\vphantom{x} } _{| _{_0} } \Big), d\gamma\Big(\frac{\partial}{\partial y}{\vphantom{y} } _{| _{_0} } \Big)\Big)=\quad
\\
\\
=\gamma^*i\partial_{_J }\bar{\partial}_{_J }f
\,\Big(\frac{\partial}{\partial x}{\vphantom{x} } _{| _{_0} } , \frac{\partial}{\partial y}{\vphantom{y} } _{| _{_0} } \Big)
=i\partial_{_j }\bar{\partial}_{_j}(f\circ\gamma)\,
\Big(\frac{\partial}{\partial x}{\vphantom{x} } _{| _{_0} } , \frac{\partial}{\partial y}{\vphantom{y} } _{| _{_0} } \Big)=
\\
\\
=i\frac{\partial^2(f\circ\gamma)}{\partial z\,\partial\bar{z} } (0)\,dz\wedge d \bar{z}\,
\Big(\frac{\partial}{\partial x}{\vphantom{x} } _{| _{_0} } , \frac{\partial}{\partial y}{\vphantom{y} } _{| _{_0} } \Big)=\frac{1}{2}\Delta (f\circ\gamma)(0) \quad 
\end{eqnarray*}
ce qui montre la nécessité de la positivité. Le même calcul, avec $z\in B_{\delta} $ à la place de $0$, montre aussi la suffisance de la positivité.\hfill $\Box$
\\
\begin{defi}
Une fonction $f\in {\cal C} ^{2}(X,\R)$ sur une variété presque complexe $(X,J)$ est dite strictement $J$-plurisousharmonique s'il existe une métrique hermitienne $\omega \in {\cal C}^0(\Lambda ^{1,1}_{_J}T_X^*)(X)$ sur le fibré tangent telle que $i\partial_{_J }\bar{\partial}_{_J }f\geq \omega $.
\end{defi}  
{\bf Quelques exemples élémentaires de fonctions strictement $J$-plurisousharmonique}. 
\\
{\bf Exemple 1}. 
On déduit facilement que si $(z_1,...,z_n)$ sont des coordonnées $\ci$ telles que $J(0)=J_0$, on a l'écriture
$$
i\partial_{_J }\bar{\partial}_{_J }f\,(\xi,J\xi)(z)=2\sum_{k,l} \frac{\partial^2f}{\partial z_k\partial\bar{z}_l}(z)\,\xi _k\bar{\xi }_l+O(|z|)(\xi ,\bar{\xi }),
$$
(voir par exemple le lemme \ref{Hess2}). 
On a alors que la fonction $f(z)=|z|^2$ est strictement $J$-plurisousharmonique sur un voisinage de l'origine des coordonnées.
\\
{\bf Exemple 2}.
Soit $F^{\lambda }_{J}:=S^{\lambda}T_{X,J}$ une puissance de Schur du fibré tangent et considérons une métrique hermitienne sur $F^{\lambda }_{J}$ telle que la courbure au sens de Griffiths soit strictement négative en un point $x$. Soit $(\sigma_k)_k \subset {\cal E}(F^{\lambda }_{J} )(U)$ un repère local
presque-holomorphe spécial au point $x\in U$. En utilisant le lemme précédent on déduit, d'après les remarques faites dans \cite{Pal}, que les fonctions $f_k:=|\sigma_k |^2_h$ sont strictement $J$-plurisousharmoniques au voisinage du point $x$.\hfill $\Box$
\\
\\
Pour réduire l'hypothèse de régularité de la fonction $f$, on a besoin de donner quelques éléments de la théorie des courants positifs sur les variétés presque complexes. Pour faire ceci on a besoin de quelques résultats et notions préliminaires que nous présentons tout de suite.
\section{Plongements par feuilles courbes $J$-holomorphes et champs de vecteurs $J$-plats sur les variétés presque complexes} 
On désigne par $B^{n}_{2r}\subset \C^n$ la boule ouverte de dimension $2n$, de rayon $2r$ et de centre l'origine.
On veut plonger dans une variété presque complexe $(X,J)$ de dimension complexe $n$ le cylindre $B^1_{\delta }\times B^{n-1}_{\delta } \subset \C^n$, avec $\delta >0$ suffisamment petit, de telle sorte que les disques $B^1_{\delta }\times p,\; p\in B^{n-1}_{\delta }$ se plongent de façon $J$-holomorphe dans $X$. De plus on veut pouvoir plonger dans toutes les ``positions possibles'' les cylindres précédents. L'existence de ces plongements est  directement liée à la notion de champ de vecteurs $J$-plat qu'on introduit ci-dessous.
\begin{defi} 
Un champ de vecteurs réel $\xi \in{\cal E}(T_X\smallsetminus 0_X)(U)$ au dessus d'un ouvert $U$ est dit $J$-plat s'il vérifie l'équation différentielle non linéaire de premier 
ordre $[\xi ,J\xi]=0$ sur l'ouvert $U$. On désigne par $P_{_J}(U,T_X)$ l'ensemble des champs de vecteurs $J$-plats au dessus de $U$.
\end{defi}
D'après le théorème de Newlander-Nirenberg on déduit que si la structure presque complexe est intégrable alors tout champ de vecteurs réel holomorphe  $\xi \in{\cal O}(T_X\smallsetminus 0_X)(U)$ au dessus d'un ouvert $U$ quelconque est $J$-plat. On a le résultat général suivant qui assure la possibilité d'effectuer des plongements du cylindre, dont les feuilles sont des courbes $J$-holomorphes, en toutes les positions possibles et l'existence locale ``en grande quantité'' des champs de vecteurs $J$-plats.
\begin{theo}\label{Jplat} 
Soit $(X,J)$ une variété presque complexe de dimension complexe $n$. Pour tout point $x_0\in X$ il existe un voisinage ouvert $U_{x_0}$ de  $x_0$  et un voisinage ouvert ${\cal B}(T_{U_{x_0}})\subset T_{U_{x_0}}$, ${\cal B}(T_{U_{x_0}})\simeq U_{x_0}\times B^n$, de la section nulle sur  $U_{x_0}$
tels que:
\\
\\
{\bf A)} Il existe une application de classe ${\cal C}^{\infty}$
$$
\Phi  : B_1^1\times {\cal B}(T_{U_{x_0}})\longrightarrow X
$$
telle que pour tout $v\in {\cal B}(T_{U_{x_0}})$ l'application $z\in B_1^1\mapsto \Phi (z,v)$ est une courbe $J$-holomorphe qui vérifie la condition $\partial _t\Phi (0,v)=v,\,z=t+is$.\\
\\
{\bf B)} Il existe une famille de plongements $(\Psi_{\alpha} :B_{\delta } ^1\times B_{\delta } ^{n-1}\longrightarrow X)_{\alpha \in I}$ de classe $\ci$  telle que pour tout $\alpha \in I$ et $z_2\in B_{\delta } ^{n-1}$ les applications 
$$
z_1\in B^1_{\delta}\mapsto\Psi_{\alpha } (z_1,z_2)
$$
sont des courbes $J$-holomorphes, 
$\Psi_{\alpha} (B_{\delta } ^1\times B_{\delta } ^{n-1})\supset U_{x_0}$
et
$$
T_{X,p}\smallsetminus 0_p =\Big\{\lambda \partial _t\Psi_{\alpha}(0,0)\,|\,\lambda \in \R\smallsetminus 0,\,\alpha \in I \Big\}
=\Big\{ \xi (p)\,|\,\xi \in P_{_J}(U_{x_0},T_X) \Big\}
$$
pour tout $p\in U_{x_0}$, $(z_1=t+is)$.
\end{theo}
Avant de passer à la preuve du théorème $\ref{Jplat}$ on a besoin de quelques préliminaires techniques. Soit $J_0$ la structure presque complexe usuelle sur $\R^{2n}$ identifié avec $\C^n$ via l'identification $z\equiv (x,y)$. Nous considérons un système de coordonnées locales centrées en $x_0\in X$ et on suppose, quitte à effectuer un changement linéaire de coordonnées, que $J(0)=J_0$. On considère aussi une boule ouverte $B^{n}_{2r}\subset \C^n$ sur laquelle l'endomorphisme $J+J_0$ est inversible et on pose par définition 
$$
q_{J}:=(J_0+J)^{-1}\cdot (J_0-J)\in \ci(End_{\R}(\R^{2n} ) )(B^{n}_{2r}).
$$
On remarque que $q_J(0)=0$. On suppose pour simplifier les notations qui suivront que $r=1$. Si $\gamma :B^1_1 \longrightarrow (B^n_2, J)$ est une courbe $J$-holomorphe, la condition de $J$-holomorphie $\partial_s\gamma =J(\gamma)\partial_t\gamma,\;z=t+i\,s$ peut être écrite de façon équivalente sous la forme
\begin{eqnarray}\label{eqJ-hol}  
\partial _{\bar{z} }\gamma+ q_{J}(\gamma)\partial_z\gamma=0
\end{eqnarray} 
où $\partial _{\bar{z} }:=\frac{1}{2}(\partial_t+J_0\partial_s)$ et $\partial _z:=\frac{1}{2}(\partial_t-J_0\partial_s)$. En effet en utilisant les identités $\partial_t=\frac{1}{2}(\partial_z+\partial_{\bar{z}})$ et $\partial _s=\frac{J_0}{2}(\partial_z-\partial_{\bar{z}})$ on peut écrire la condition  $\partial_s\gamma =J(\gamma)\partial_t\gamma$ sous la forme
$$
(J_0+J(\gamma ))\,\partial_{\bar{z}}\gamma =(J_0-J(\gamma ))\,\partial_z\gamma.
$$ 
L'inversibilité de l'endomorphisme $J+J_0$ donne alors l'écriture sous la forme $\eqref{eqJ-hol}$.
On rappelle aussi (voir l'article de Sikorav dans l'ouvrage de Audin et Lafontaine \cite{Ad} pour plus de détails) que l'opérateur 
$$
P:{\cal C}^{k+\mu }(B_1^1;\C^n) \longrightarrow {\cal C}^{k+\mu +1}(B_1^1;\C^n),
$$
$\;k\in\N,\;\mu \in (0,1)$ défini par la formule $P\gamma(z):=P'\gamma(z)-P'\gamma(0)$, avec
\begin{eqnarray*}
P'\gamma(z):=\frac{1}{2\pi i}  \int\limits_{\zeta \in B_1^1}\frac{\gamma(\zeta )}{\zeta -z}\;d\zeta \wedge d\bar{\zeta},   
\end{eqnarray*}
vérifie les propriétés suivantes: $\partial _{\bar{z} }\circ P=\I$ et pour tout entier $k\in\N$ et $\mu \in (0,1)$ il existe une constante $c_{k,\mu }>0$ telle que pour 
toute courbe $\gamma \in {\cal C}^{k+\mu }(B_1^1;\C^n)$ on a l'estimation 
\begin{eqnarray}\label{Cauchy-Oper} 
\|P\gamma\|_{k+\mu +1}\leq (2+c_{k,\mu })\|\gamma\|_{k+\mu },
\end{eqnarray}
où $\|\cdot\|_{k+\mu }$ désigne la norme de Hölder usuelle sur $B^1_1$. Pour prouver le théorème $\ref{Jplat} $ on utilisera la remarque essentielle suivante, utilisée par McDuff (voir le lemme 1.4 dans \cite{McD} ) et aussi par Sikorav pour prouver le théorème 3.1.1 dans l'ouvrage \cite{Ad}: une courbe  $\gamma :B^1_1 \longrightarrow (B^n_2, J)$ est $J$-holomorphe si et seulement si la courbe 
$$
\gamma _0:=\gamma +P\Big(q_J(\gamma)\,\partial_z\gamma\Big)
$$
 est $J_0$-holomorphe. De plus on a l'égalité $\gamma _0(0)=\gamma(0)$.
\\ 
\\
On aura besoin de quelques remarques élémentaires de topologie différentielle qui seront utilisées plusieurs fois dans la suite.
\\
\\
$Remarque\; 1$. Soit $f:X\times Y\longrightarrow Z$ une application entre espaces topologiques telle que l'application 
$\Phi_f:X\times Y\longrightarrow X\times Z,\;\Phi_f(x,y):=(x,f(x,y))$ soit ouverte. Alors pour tout $(x_0,y_0)\in X\times Y $, pour tout voisinage ouvert $V_{y_0}\subset Y$ de $y_0$ et pour tout compact $K\subset f_{x_0}(V_{y_0})$ (ici on pose par définition $f_x:=f(x,\cdot)$) il existe un voisinage ouvert $U_{x_0}\subset X$ de $x_0$ tel que
pour tout $x\in U_{x_0}$ on a   $f_x(V_{y_0})\supset K$. L'hypothèse précédente est vérifiée par exemple si $f$ est une application de classe ${\cal C}^1$ entre variétés de Banach telle que pour tout $x\in X$ l'application $f_x:Y\longrightarrow Z$ soit un plongement ouvert, autrement dit $f_x$ est injective et 
$$
d_yf_x:T_{Y,y} \longrightarrow T_{Z,f_x(y)}  
$$
est un isomorphisme pour tout $y\in Y$. En effet dans ce cas le théorème d'inversion locale implique que l'application $\Phi_f$ est ouverte. 
\\
\\
$Remarque\; 2$. 
Dans le cas où l'application $f_{x_0} :Y\longrightarrow Z$ est un plongement ouvert seulement en un point $x_0\in X$ on a d'après le théorème des fonctions implicites que pour tout compact $K\subset Z$ il existe un voisinage ouvert $V_K\subset Z$ de $K$, un voisinage ouvert $W\subset Y$ de $f^{-1}_{x_0} (V_K)$ et un voisinage ouvert $U_{x_0}\subset X$ de $x_0$ tel que pour tout $x\in U_{x_0}$ l'application 
$$
f_x:f^{-1}_x(V_K)\cap W \longrightarrow V_K
$$ 
soit un difféomorphisme de classe ${\cal C}^1$.
\\
\\
$Remarque\; 3$. Le théorème des fonctions implicites implique que si $f:X\times Y'\longrightarrow Z$ est une application de classe ${\cal C}^1$ entre variétés de Banach telle qu'il existe un point $x_0\in X$ et un ouvert relativement compact $Y\subset Y'$ (donc $Y'$ est de dimension finie) tels que $f_{x_0}:Y\longrightarrow Z $ soit injective et 
$$
d_yf_{x_0} :T_{Y',y} \longrightarrow T_{Z,f_{x_0} (y)}
$$
soit un isomorphisme pour tout $y\in \overline{Y}$, alors il existe un voisinage ouvert $U_{x_0}\subset X$ de $x_0$ tel que
pour tout $x\in U_{x_0}$ l'application $f_x:Y\longrightarrow Z$ est un plongement ouvert.
\\
\\
{\bf Preuve du théorème \ref{Jplat}} 
\\
{\bf Preuve de la partie A} 
\\
Pour tout entier $k\in \N,\;k\geq 2$ nous considérons l'application de classe ${\cal C}^{k-1} $
\begin{diagram}[height=1cm,width=1cm]
&F:[0,1]\times {\cal C}^{k+\mu }&(B_1^1\times B_1^n\times B_1^n;\,B_2^n)&\rTo&{\cal C}^{k+\mu }(B_1^1\times B_1^n\times B_1^n;\,\C^n)
\\
&&(\varepsilon ,\phi)&\longmapsto& \phi+P_{z}\Big(q_J(\varepsilon \phi)\,\partial_{z}\phi\Big)
\end{diagram}
où $\mu \in (0,1)$ est une constante fixée et $(P_{z}\phi)(z,x,v):=(P\phi(\cdot,x,v))(z),\;(z,x,v)\in B_1^1\times B_1^n\times B_1^n$. Considérons aussi l'application holomorphe
$H\in {\cal O} (B_1^1\times B_1^n\times B_1^n;\,B_2^n)$ définie par la formule $H(z,x,v):=x+zv$.
\\
Le fait que l'application 
$$
F_0:=F(0,\cdot):{\cal C}^{k+\mu }(B_1^1\times B_1^n\times B_1^n;B_2^n)\longrightarrow {\cal C}^{k+\mu }(B_1^1\times B_1^n\times B_1^n;\C^n)
$$
soit l'inclusion canonique entraîne, d'après la remarque 2, l'existence d'un voisinage ouvert 
${\cal V}_k\subset {\cal C}^{k+\mu }(B_1^1\times B_1^n\times B_1^n;B_2^n)$ de $H$ (avec
${\cal V}_k\supset {\cal V}_{k+1} $), d'un voisinage ouvert 
$$
{\cal W}_k\subset {\cal C}^{k+\mu }(B_1^1\times B_1^n\times B_1^n;B_2^n)
$$
de ${\cal V}_k$, ${\cal W}_k\supset {\cal W}_{k+1}$  et de $\varepsilon _0\in (0,1]$ tel que pour tout $\varepsilon \in [0,\varepsilon _0]$ l'application 
$$
F_{\varepsilon}:F_{\varepsilon}^{-1} ({\cal V}_k)\cap {\cal W}_k \longrightarrow {\cal V}_k
$$
est un difféomorphisme de classe ${\cal C}^{k-1}$. On pose alors par définition 
$\phi_{\varepsilon}:=F^{-1}_{\varepsilon}(H)$ et on remarque que l'application 
$
\varepsilon \,\phi_{\varepsilon}\in \ci(B_1^1\times B_1^n\times B_1^n;B_2^n)
$ est $J$-holomorphe par rapport à la variable $z\in B_1^1$.  Nous considérons maintenant l'application de classe $\ci$
\begin{diagram}[height=1cm,width=1cm]
&\chi:[0,&\varepsilon _0]\times B_1^n\times B_1^n&\rTo& B_2^n\times \C^n&&
\\
& &(\varepsilon ,x,v)&\longmapsto&(x,\partial_t\phi_{\varepsilon}(0,x,v)), & 
\end{diagram}
$z=t+i\,s$. On rappelle que $\phi_{\varepsilon}(0,x,v)=x$. Le fait que l'application $\chi_0:=\chi(0,\cdot,\cdot)$ soit l'inclusion canonique entraîne, d'après les remarques 3 et 1, que quelque soit $r\in (0,1)$ et
$\rho \in (0,r)$ il existe $\varepsilon _1=\varepsilon _1(r,\rho )\in(0,\varepsilon _0]$ tel que pour tout $\varepsilon \in [0,\varepsilon _1]$ l'application 
$$
\chi_{\varepsilon }:B_r^n\times B_r^n\longrightarrow \chi_{\varepsilon }(B_r^n\times B_r^n)\supset \overline{B_{\rho}^n}\times \overline{B_{\rho}^n}
$$
est un difféomorphisme de classe $\ci$. On considère l'application 
$\chi_{\varepsilon }^{-1}:\overline{B_{\rho}^n}\times \overline{B_{\rho}^n}\longrightarrow  B_r^n\times B_r^n $ et on définit l'application
\begin{diagram}[height=1cm,width=1cm]
&\Phi_{\varepsilon} :& B_1^1\times B_{\varepsilon \rho}^n \times B_{\varepsilon \rho}^n
 &\rTo& B_2^n&
\\
& &(z,x,v)&\longmapsto&\varepsilon \phi_{\varepsilon}(z,\chi_{\varepsilon }^{-1}(\varepsilon^{-1} x,\varepsilon^{-1} v)).&
\end{diagram}
Si on pose par définition $U_{x_0}:=B_{\varepsilon \rho}^n$ et ${\cal B}(T_{U_{x_0}}):=B_{\varepsilon \rho}^n \times B_{\varepsilon \rho}^n$ on a que l'application $\Phi_{\varepsilon}$ vérifie les conditions de la partie A de l'enoncé du théorème \ref{Jplat}. On verra de suite que pour satisfaire aussi la conclusion B du théorème \ref{Jplat} il est nécessaire de considérer un voisinage ouvert $U_{x_0}$ plus petit.
\\
\\
{\bf Preuve de la partie B} 
\\
On rappelle qu'on désigne par $J_0$ la structure presque complexe usuelle sur $\R^{2n}$ identifié à $\C^n$ via $z\equiv (x,y)$. Avec cette identification on voit le groupe $Gl(n,\C)$ comme sous-groupe du groupe $Gl(2n,\R)$. Précisément 
$A\in Gl(n,\C)\subset Gl(2n,\R)$ si et seulement si $A\,J_0=J_0\,A$. Dans la suite on désignera par $U(n):=O(2n)\cap Gl(n,\C)$ le groupe unitaire. Soit $\delta =\rho /2$ et 
$$
l:B_{\delta }^n\times B_{\delta}^{n-1}\times U(n)\longrightarrow B_{\rho }^n \subset \R^{2n}\equiv \C^n    
$$
l'application définie par la formule
\begin{eqnarray*}
l(p,z_2,A)=p+A\Big(x_2\cdot \frac{\partial}{\partial x_2}+y_2\cdot\frac{\partial}{\partial y_2}\Big),
\end{eqnarray*}
avec l'identification $z_2\equiv(x_2,y_2)$ et $x_2\equiv(x_2,...,x_n),\;y_2\equiv(y_2,...,y_n)$. 
Considérons donc l'application de classe $\ci$
\begin{diagram}[height=1cm,width=1cm]
&\Phi :&[0,\varepsilon _1]\times (B_{\delta } ^1\times B_{\delta } ^{n-1})\times B_{\delta }^n\times U(n)  &\rTo& B_2^n&
\\
& &(\varepsilon ;(z_1,z_2);(p,A))&\longmapsto&\phi_{\varepsilon} \Big(z_1, 
\chi _{\varepsilon}^{-1}\Big(l(p,z_2,A);A\frac{\partial}{\partial x_1}\Big)\Big), &
\end{diagram}
où $\phi_{\varepsilon}:=F^{-1}_{\varepsilon}(H)$ est l'application définie dans la preuve de la partie A. Avec l'identification $\Phi(\varepsilon ;(z_1,z_2);(p,A))\equiv \Phi ^{p,A}_{\varepsilon }(z_1,z_2)$ on a les propriétés suivantes:
$$
\left  \{
\begin{array}{lr}
\Phi ^{p,A}_{\varepsilon }(0,z_2)=l(p,z_2,A)
\\
\\
\partial_t \Phi ^{p,A}_{\varepsilon }(0,z_2)=\displaystyle { A\frac{\partial}{\partial x_1}}
\\
\\
\Phi ^{p,A}_0(z_1,z_2)=\displaystyle {t\,A\frac{\partial}{\partial x_1}+s\,A\frac{\partial}{\partial y_1}+l(p,z_2,A)}
\end{array}
\right.
$$
(rappelons que $z_1:=t+i\,s$). Soit $\delta _1\in (0,\delta )$ un réel suffisamment petit pour pouvoir assurer l'inclusion 
$\overline{B^n_{\delta _1} }\subset \Phi ^{p,A}_0(B_{\delta } ^1\times B_{\delta } ^{n-1})$ pour tout $p\in\overline{B^n_{\delta _1} }$ et  $A\in U(n)$. On aura alors que l'image du plongement 
\begin{diagram}[height=1cm,width=1cm]
&\Phi_0\times \I :& (B_{\delta } ^1\times B_{\delta } ^{n-1})\times\overline{B^n_{\delta _1} }\times U(n)  &\rTo& B_2^n\times\overline{B^n_{\delta _1} }\times U(n)&
\\
& &((z_1,z_2);(p,A))&\longmapsto&(\Phi ^{p,A}_0(z_1,z_2);(p,A)) &
\end{diagram}
contient le compact $\overline{B^n_{\delta _1} }\times\overline{B^n_{\delta _1} }\times U(n)$ (rappelons que le groupe $U(n)$ est compact). On aura d'après les remarques 3 et 1 l'existence de $\varepsilon _2\in (0,\varepsilon _1]$ tel que pour tout $\varepsilon \in (0,\varepsilon _2]$ et $(p,A)\in B^n_{\delta _1}\times U(n)$ l'application
$$
\Phi ^{p,A}_{\varepsilon } :B_{\delta } ^1\times B_{\delta } ^{n-1}\longrightarrow \Phi^{p,A}_{\varepsilon }(B_{\delta } ^1\times B_{\delta } ^{n-1})\supset B^n_{\delta _1}
$$
est un difféomorphisme de classe $\ci$. Nous considérons donc l'application 
$$
\Psi_{\varepsilon}: (B_{\delta } ^1\times B_{\delta } ^{n-1})\times B^n_{\varepsilon \delta _1}\times U(n)\longrightarrow B_2^n
$$
définie par la formule $\Psi^{p,A}_{\varepsilon } :=\varepsilon \,\Phi ^{\varepsilon ^{-1} p,A}_{\varepsilon }$, $(p,A)\in B^n_{\varepsilon \delta _1}\times U(n)$ et on remarque qu'elle vérifie les propriétés suivantes:
\begin{eqnarray*}
\partial_{\bar{z}_1} \Psi^{p,A}_{\varepsilon }+q_J(\Psi^{p,A}_{\varepsilon })\,\partial_{z_1} \Psi^{p,A}_{\varepsilon }=0&\qquad&\Psi^{p,A}_{\varepsilon }(0,0)=p
\\
\\
\Psi^{p,A}_{\varepsilon }(B_{\delta } ^1\times B_{\delta } ^{n-1})\supset B^n_{\varepsilon \delta _1}&\qquad&\partial_t \Psi^{p,A}_{\varepsilon }(0,z_2)=
\displaystyle { \varepsilon A\frac{\partial}{\partial x_1}}
\end{eqnarray*}
On définit les champs de vecteurs $J$-plats
\begin{eqnarray*}
\xi _{p,A}:=\partial_t\Psi^{p,A}_{\varepsilon }\circ (\Psi^{p,A}_{\varepsilon })^{-1}  &\qquad&\xi _{p,A}(p)=\displaystyle { \varepsilon A\frac{\partial}{\partial x_1}}
\end{eqnarray*}
 sur l'ouvert $B^n_{\varepsilon \delta _1}$. Le fait que l'action de $U(n)$ est transitive sur la sphère $S^{2n-1}$, (voir \cite{Bo-Tu}) entraîne que la famille 
$$
\{\lambda \xi _{p,A} \in P_J(B^n_{\varepsilon \delta _1}, T_X )\;|\; \lambda \in \R\smallsetminus \{0\},\, (p,A)\in B^n_{\varepsilon \delta _1}\times U(n) \} 
$$ 
engendre ponctuellement (au sens ensembliste) $T_{X|B^n_{\varepsilon \delta _1} }\smallsetminus0_X$, ce qui prouve la partie B du théorème $\ref{Jplat}$ avec $U_{x_0}:=B^n_{\varepsilon \delta _1} $ et $I:=B^n_{\varepsilon \delta _1}\times U(n)$.\hfill $\Box$
\\
\\
Le lemme élémentaire suivant montre que tout champ de vecteurs $J$-plat provient localement d'un plongement dont les feuilles  sont des courbes $J$-holomorphes.
\begin{lem}\label{locJ-plat}Soit $(X,J)$ une variété presque complexe de dimension complexe $n$ et $\xi$ un champ de vecteurs $J$-plat sur un ouvert $U$. Pour tout $x\in U$ il existe un voisinage ouvert $U_x\subset U$ de $x$ et une carte locale 
$(U_x,\sigma _{\xi }^{-1}),\,\sigma_{\xi}  :B^1_{\delta}\times B^{n-1}_{\delta  }\longrightarrow U_x$, compatible avec l'orientation canonique de $(U_x, J)$ telle que pour tout 
$z_2\in B^{n-1}_{\delta}$, les applications 
$z_1\in B^1_{\delta}\mapsto \sigma_{\xi}(z_1,z_2)$ sont des courbes $J$-holomorphes et 
$d\sigma _{\xi } ( \frac{\partial}{\partial x_1} )=\xi \circ\sigma _{\xi },\,z_1=x_1+iy_1$.
\end{lem}
$Preuve$. Soient $v_2,...,v_n \in T_{X,x}$ des vecteurs tels que $\xi (x),v_2,...,v_n$ soit une base sur $\C$ de $T_{X,x}$ et $\tau$ des coordonnées locales centrées en $x$ telles que:
\begin{eqnarray*} 
d\tau^{-1}\Big(\frac{\partial}{\partial x_1}{\vphantom{x_1} } _{| _{_0} } \Big)=\xi (x) && d\tau^{-1}\Big(\frac{\partial}{\partial y_1}{\vphantom{y_1} } _{| _{_0} } \Big)=J\xi (x)
\\
\\
d\tau^{-1}\Big( \frac{\partial}{\partial x_k}{\vphantom{x_k} } _{| _{_0} } \Big)=v_k && d\tau^{-1}\Big( \frac{\partial}{\partial y_k}{\vphantom{y_k} } _{| _{_0} } \Big)=Jv_k
\end{eqnarray*}
pour tout $k=2,...,n$, (on désigne par $(x_1,y_1,...,x_n,y_n)$ les coordonnées sur $\R^{2n}$). On désigne par $\phi_{\xi},\,\phi_{J\xi}:V_x\times (-\delta ,\delta )\subset X\times \R\longrightarrow X $ les flots respectifs des champs $\xi $ et $J\xi$ au voisinage $V_x$ de $x$ (pour simplifier les notations on utilisera dans la suite l'identification $\phi_{\xi}(x,t)\equiv\phi^t_{\xi}(x) $) et on considère l'application $\sigma _{\xi }:\mbox{Im}\,\tau\longrightarrow X$ définie par la formule
\begin{eqnarray*} 
\sigma _{\xi }(x_1,y_1,...,x_n,y_n)&:=&\phi^{x_1} _{\xi}\circ \phi^{y_1} _{J\xi}\circ\tau^{-1}(0,0,x_2,y_2,...,x_n,y_n)=
\\
\\
&=&\phi^{y_1}_{J\xi}\circ \phi^{x_1}_{\xi}\circ\tau^{-1}(0,0,x_2,y_2,...,x_n,y_n).
\end{eqnarray*}
D'après le théorème d'inversion locale on a l'existence d'un voisinage ouvert $U _x\subset U$ de $x$ tel que $(U_x, \sigma _{\xi }^{-1} )$
 soit une carte locale compatible avec l'orientation canonique de $(U_x,J)$ telle que 
\begin{eqnarray*} 
d\sigma _{\xi }  \Big( \frac{\partial}{\partial x_1} \Big)=\xi \circ\sigma _{\xi } \quad \text{et}\quad 
d\sigma _{\xi }  \Big( \frac{\partial}{\partial y_1} \Big)=J\xi \circ\sigma _{\xi } .
\end{eqnarray*}
Si on suppose $\sigma _{\xi }^{-1}(U_x)=B^1_{\delta}\times B^{n-1}_{\delta  } \subset \R^2\times \R^{2n}$ on en déduit que les applications 
$(t,s)\in B^1_{\delta}\mapsto \sigma _{\xi }(t,s,a_2,b_2,...,a_n,b_n)$ sont des courbes $J$-holomorphes pour tout $(a_2,b_2,...,a_n,b_n)\in B^{n-1}_{\delta  }$.\hfill $\Box$
\\
\\
On aura besoin aussi du lemme suivant.
\begin{lem}\label{perturbCurbJ-hol}  
Soit $(X,J)$ une variété presque complexe de dimension complexe $n$ et soit $\gamma :B^1_{\delta }\longrightarrow X$ une courbe $J$-holomorphe lisse. Il existe alors un plongement 
\\
$\sigma :B^1_{\rho }\times B^{n-1}_{\rho }\longrightarrow X,\,\rho \in (0,\delta )$ de classe $\ci$ qui préserve les orientations canoniques  tel que les applications $\sigma (\cdot ,z_2),\,z_2\in B^{n-1}_{\rho }$ soient des courbes $J$-holomorphes et $\sigma (\cdot , 0)=\gamma$.
\end{lem} 
$Preuve$. Soit $B^n_2\subset X$ une boule coordonnée telle que $J(0)=J_0$ et $\gamma (0)=0$. Soit $\mu _{\lambda }:B^1_1\longrightarrow B^1_{\lambda  },\,\mu _{\lambda }(z)=\lambda z$ l'homothétie de facteur $\lambda >0$ et $\gamma _{\lambda },\,\lambda \in (0,\delta]$ la courbe $J$-holomorphe définie par la formule $\gamma _{\lambda }:=\gamma \circ \mu _{\lambda }$. Considérons la famille de courbes $J_0$-holomorphes $(u_{\lambda})_{\lambda \in (0,\delta]}$ définie par la formule
$$
u_{\lambda}:=\lambda ^{-1}\Big[\gamma _{\lambda }+P\Big(q_{J}(\gamma _{\lambda } )\,\partial_z  \gamma _{\lambda }\Big)\Big]. 
$$
Considérons des vecteurs $\xi _2,...,\xi _n\in \R^{2n}\equiv\C^n$ tels que les vecteurs 
$\lambda \partial_t u_{\lambda}(0),\xi _2,...,\xi _n$ forment une base $J_0$-complexe de $\R^{2n}$ et la famille d'applications $J_0$-holomorphes 
$$
(H_{\lambda } )_{\lambda \in (0,\delta ]}\subset {\cal O} (B^1_1\times B^{n-1}_1;B^n_2),
$$
définie par la formule $H_{\lambda }(z_1,z_2)=u_{\lambda}(z_1)+\xi \cdot z_2$. Nous considérons aussi l'application
$$
F:[0,1]\times {\cal C}^{k+\mu }(B_1 ^1\times B_1 ^{n-1} ;\,B_2^n)\longrightarrow {\cal C}^{k+\mu }(B_1 ^1\times B_1 ^{n-1} ;\,\C^n)
$$
définie comme dans la preuve du théorème $\ref{Jplat}$. Le fait que l'ensemble 
$$
\overline{(H_{\lambda } )_{\lambda \in (0,\delta ]}}\subset {\cal C}^{k+\mu }(B_1 ^1\times B_1 ^{n-1} ;\,B_2^n) 
$$ 
soit compact, (pour tout $k\geq 1$) entraîne, d'après la remarque 2 de la preuve du théorème $\ref{Jplat}$, l'existence d'un $\rho \in (0,\delta]$ pour lequel il existe les applications $\phi _{\varepsilon }:= F^{-1}_{\varepsilon }(H_{\varepsilon }) ,\,\varepsilon \in (0,\rho ]$, (les applications 
$F^{-1}_{\varepsilon }$ sont  définies comme dans la preuve du théorème $\ref{Jplat}$). De façon explicite on a donc l'identité
$$
\phi_{\varepsilon}+P_{z_1}\Big(q_J(\varepsilon \phi_{\varepsilon} )\,\partial_{z_1}\phi_{\varepsilon}\Big) =H_{\varepsilon }.
$$ 
On déduit alors, grâce à l'inégalité $\eqref{Cauchy-Oper}$, que pour $\varepsilon >0$ suffisamment petit, (disons $\varepsilon \in (0,\rho ]$), on a l'inégalité
$$
\|\phi_{\varepsilon}-H_{\varepsilon }\|_{k+\mu +1}\leq\varepsilon C_{k,\mu }\|d_0q_{J} \|\cdot \|\phi_{\varepsilon} \|_{k+\mu}\cdot \|d\phi_{\varepsilon} \|_{k+\mu},     
$$
qui compte tenu de la compacité de la famille 
$
\overline{(\phi_{\varepsilon })_{\varepsilon \in (0,\rho] }}
\subset {\cal C}^{k+\mu }(B_1 ^1\times B_1 ^{n-1} ;\,B_2^n) 
$, (pour tout $k\geq 1$) implique l'inégalité
$$
\|\phi_{\varepsilon}-H_{\varepsilon }\|_{k+\mu +1}\leq C'_{k,\mu }\varepsilon.
$$
On considère le plongement linéaire $L(z_1,z_2):=d_0\gamma (z_1)+\xi \cdot z_2$ et on remarque l'inégalité 
$$
\|H_{\varepsilon }-L\|_{k+\mu+1 }=\|u_{\varepsilon}-d_0\gamma \|_{k+\mu+1 }\leq \varepsilon C'_{k,\mu }\|d_0q_{J}\|\cdot \|d_0\gamma\|^2,
$$
pour tout  $\varepsilon \in (0,\rho ]$. On déduit alors que les applications $\phi_{\varepsilon }$ sont des plongements pour $\rho >0$ suffisamment petit (voir lemme 1.3 du chapitre 2 dans l'ouvrage de Hirsch \cite{Hir}) . On considère donc les plongements $\psi_{\varepsilon }:=\varepsilon \phi_{\varepsilon }$ et on remarque les égalités
$$
\Big[\psi_{\varepsilon}+P_{z_1}\Big(q_J(\psi_{\varepsilon} )\,\partial_{z_1}\psi_{\varepsilon}\Big) \Big](\cdot ,0)=
\varepsilon u_{\varepsilon }=
\gamma _{\varepsilon}+P_{z_1}\Big(q_J( \gamma _{\varepsilon} )\,\partial_{z_1}\gamma _{\varepsilon}\Big) 
$$
qui montrent l'égalité $\psi_{\varepsilon }(\cdot,0)=\gamma _{\varepsilon}$. On déduit alors que l'application
$$
(z_1,z_2)\in B_{\rho  } ^1\times B_{\rho } ^{n-1}\mapsto 
\sigma (z_1,z_2):=\psi_{\rho  }(\rho^{-1} z_1,z_2)
$$
est le plongement  voulu.\hfill $\Box$
\section{Courants positifs sur les variétés presque complexes}
\subsection{Généralités} 
On commence par rappeler quelques définitions générales de la théorie des courants.
\begin{defi}\label{mass}  
Soit $\Theta \in {{\cal D}'}^k(X)$ un courant de degré $k$, d'ordre zéro sur une variété différentiable $X$ orientable et orientée de dimension $n$. Une masse du courant $\Theta $ est une mesure de Radon positive $\mu$ sur $X$ telle que si $\psi\in{\cal E}(\Lambda ^n T_X^*)(X)$ est une forme de volume arbitraire et si $A\subset X$ est un ensemble de Borel alors $\mu (A)=0$ si et seulement si $\int _A \Theta(\xi _1,...,\xi _k)\cdot \psi=0$  pour tout champ de vecteurs $\xi _1,...,\xi _k\in {\cal E}(T_X)(X)$.
\end{defi}
On remarque que si $\mu _1$ et $\mu _2$ sont deux masses du même courant alors l'une est absolument continue par rapport à l'autre. Il est bien connu, (cf. \cite{Fed}, \cite{G-M-S}) que tout courant d'ordre zéro admet une masse qui peut être définie par la formule
\begin{eqnarray*}
\mu _g(\Theta )(U):=\sup_{\displaystyle
\scriptstyle   \varphi \in {\cal D}^{n-k} (U)  
\atop
\scriptstyle |\varphi|_g \leq 1 } \left|  \,\int\limits_{U}\Theta \wedge \varphi  \right| 
\end{eqnarray*}
pour tout ouvert $U\subset X$ relativement compact dans $X$, (ici $g$ est une métrique Riemannienne sur $X$). Avec les notations de la définition $\ref{mass}$ on a par conséquence du Théorème de Radon-Nikodym l'existence d'une $k$-forme $\theta _{\mu ,\psi} $ telle que  pour tout champ de vecteurs  $\xi _1,...,\xi _k\in {\cal E}(T_X)(X)$ la fonction 
$\theta_{\mu ,\psi}(\xi _1,...,\xi _k)\in L^1_{loc}(X,{\cal B}_X,\mu )$ (ici ${\cal B}_X$ désigne la $\sigma $-algèbre de Borel) est définie $\mu$-presque partout par la formule
\begin{eqnarray*} 
\theta_{\mu ,\psi}(\xi _1,...,\xi _k)(x):=\lim_{r\rightarrow 0 }\,\frac{1}{\mu (B_r(x))}\int\limits_{B_r(x) }\Theta(\xi _1,...,\xi _k)\cdot \psi
\end{eqnarray*}
où $B_r(x)$ est une boule de rayon $r$ relative à un ouvert coordonné quelconque. On aura alors pour tout Borelien $A\in{\cal B}_X$ l'égalité
\begin{eqnarray*} 
\int\limits_{A}\Theta(\xi _1,...,\xi _k)\cdot \psi=\int\limits_{A} \theta_{\mu ,\psi} (\xi _1,...,\xi _k)d\mu 
\end{eqnarray*}
qu'on dénote souvent sous la forme $\Theta =\theta_{\mu ,\psi}\cdot \mu $. 
Nous rappelons maintenant quelques résultats de base de la théorie des courants d'ordre zéro (cf. \cite{Fed}, \cite{G-M-S}).
\begin{theo}{\bf (Compacité faible de la masse)}. Soit $ \{\Theta _{\nu}  \}_{\nu } \subset {{\cal D}'}^k(X)$ une suite de courants d'ordre zéro telle que $\sup_{\nu }\mu (\Theta _{\nu } )(U)<\infty$ pour tout ouvert relativement compact $U$ de $X$. Il existe alors une sous-suite $\{\Theta _{\nu_j}  \}_{\nu_j }$ de $\{\Theta _{\nu}  \}_{\nu }$ convergente pour la topologie faible des courants d'ordre zéro vers un courant d'ordre zéro $\Theta   \in {{\cal D}'}^k(X)$.
\end{theo} 
Ce théorème est juste une conséquence du théorème classique de Banach-Alaoglu. Le théorème précédent admet un réciproque que nous énonçons sous la forme suivante.
\begin{theo}\label{invcompmas} 
Soit $ \{\Theta _{\nu}  \}_{\nu } \subset  {{\cal D}'}^k(X)$ une suite de courants d'ordre zéro telle que $\sup_{\nu }|\left<\Theta _{\nu},\varphi \right>|<\infty $ pour toute forme à support compact $\varphi \in {\cal C}^0(\Lambda ^{n-k}T_X^*)(X)$. Alors les masses des courants $\Theta _{\nu}$ sont localement équi-bornées au sens suivant : pour tout ouvert relativement compact $U$ de $X$on a $\sup_{\nu }\mu (\Theta_{\nu } )(U)<\infty$.
\end{theo}
Ce théorème est simplement une conséquence du théorème classique de Banach-Steinhaus. Nous avons aussi la lemme très utile suivant.
\begin{lem}\label{critconvzero} 
Soit $\{\Theta _{\nu}  \}_{\nu } \subset {{\cal D}'}^k(X)$ une suite de courants d'ordre zéro convergente faiblement vers un courant d'ordre zéro $\Theta\in {{\cal D}'}^k(X)$. Si $\sup_{\nu }\mu (\Theta_{\nu } )(X)<\infty$ alors la suite $\{\Theta _{\nu}  \}_{\nu }$ converge vers le courant $\Theta$ dans la topologie faible des courants d'ordre zéro.
\end{lem}  
Considérons à partir de maintenant  une variété presque complexe $(X,J)$ de classe ${\cal C}^{\infty}$ et de dimension réelle $2n$ munie d'une métrique $\omega \in {\cal E}(\Lambda ^{1,1}_{_J}T_X^*)(X)$ et les $(p,p)$-formes fortement positives 
$\omega_p:=1/p!\,\omega^p$ pour $p=0,...,n$. On remarque que $\omega_n$ est la forme de volume associée à la métrique $\omega$. Les notations précédentes seront utiles pour montrer l'équivalence des définitions suivantes.
\begin{defi}\label{corrpos} 
Un courant $\Theta \in {\cal D}'_{p,p}(X)$ sur une variété presque complexe $(X,J)$ est dit positif si il vérifie une des trois propriétés équivalentes suivantes.
\\
$a)$ Pour tout champ de vecteurs réels $\xi _1,...,\xi _{n-p} \in {\cal E}(T_X)(X)$ et pour toute forme $\varphi \in {\cal D}^{2n}(X)$ positive on a l'inégalité
$$
\langle\Theta(\xi _1,J\xi _1,...,\xi _{n-p} ,J\xi _{n-p}  ),\varphi  \rangle\geq 0 .
$$
\\
$b)$ Le courant $\Theta$ est d'ordre zéro et le courant $\Theta \wedge \omega _p$ détermine une masse $\|\Theta \|_{\omega}$ du courant $\Theta$ telle que quel que soit le représentant 
$$\theta _{\omega}\in \theta _{_{\|\Theta \|_{\omega}, \omega _n} }\in 
\Big({\cal E}(\Lambda ^{n-p,n-p}_{_J}\,T_X^*) \otimes_{_{{\cal E}_X(\C) }}{\cal L}^1_{loc}({\cal B}_X, \|\Theta \|_{\omega} )\Big)(X)
$$
de la forme $\theta _{_{\|\Theta \|_{\omega}, \omega _n}}$ on a que la forme $\theta _{\omega }(x)\in \Lambda ^{n-p,n-p}_{_J}\,T_{X,x} ^*$ est positive pour $\|\Theta \|_{\omega}$-presque tout $x\in X$.
\\
$c)$ Pour tout $(p,p)$-forme $\varphi  \in {\cal E}(\Lambda ^{p,p}_{_J}T_X^*)(X)$ fortement positive le courant $\Theta \wedge \varphi $ détermine une mesure de Radon positive.
\\
Le cône des courants positifs de bidimension $(p,p)$ sera noté par ${\cal D}'_{p,p}(X)^+$.
\end{defi}
$Preuve \; de\; l'\acute{e}quivalence$. Nous montrons les implications $a)\Longrightarrow c)$ et $c)\Longrightarrow b)$. L'implication $b)\Longrightarrow a)$ est évidente. Commençons par prouver l'implication $a)\Longrightarrow c)$.
\\
Soit $U\subset X$ un ouvert coordonée, soit 
$(\zeta _1,...,\zeta _n)$ un repère du fibré $ T^{1,0}_{X,J\,|U}$ et $(\rho _{\varepsilon})_{\varepsilon >0}$ une famille de noyaux régularisants usuels. Si
$$
\Theta =i^{(n-p)^2}\sum_{|K|=|H|=n-p}\Theta _{K,H}\,\zeta _K^*\wedge\bar{\zeta}_H^*    
$$
est l'expression locale du courant $\Theta$ on définit les $(n-p,n-p)$-formes
$$
\Theta*\rho _{\varepsilon} :=i^{(n-p)^2}\sum_{|K|=|H|=n-p}\Theta _{K,H}*\rho _{\varepsilon}\,\zeta _K^*\wedge\bar{\zeta}_H^*.    
$$
Soient de plus $\xi _1,...,\xi _{n-p} \in {\cal E}(T_X)(U)$ des champs de vecteurs à coefficients constants par rapport au repère $(\zeta _1,...,\zeta _n)$. L'égalité  
$$
\Theta(\xi _1,J\xi _1,...,\xi _{n-p} ,J\xi _{n-p}  )*\rho _{\varepsilon}=(\Theta*\rho _{\varepsilon} )(\xi _1,J\xi _1,...,\xi _{n-p} ,J\xi _{n-p})
$$
entraîne que les formes $\Theta*\rho _{\varepsilon}$ sont positives. Pour tout $(p,p)$-forme $\varphi  \in {\cal E}(\Lambda ^{p,p}_{_J}T_X^*)(X)$ fortement positive on a alors l'inégalité $(\Theta *\rho _{\varepsilon})\wedge \varphi \geq 0$. En passant à la limite on obtient la conclusion voulue. 
\\
\\
Nous montrons maintenant l'implication $c)\Longrightarrow b)$.
Montrons d'abord que le courant $\Theta $ est d'ordre zéro. Soit $U\subset X$  sur lequel $T_{X|U}$ est trivial et soit $(\zeta _1,...,\zeta _n)$  un repère $\omega$-orthonormé du fibré $ T^{1,0}_{X,J\,|U}$. Les formes $\omega _p$ s'expriment alors par rapport au repère choisi sous la forme
$$
\omega _p =\frac{i^{p^2}}{2^p} \sum_{|K|=p}\zeta _K^*\wedge\bar{\zeta}_K^*.    
$$
Pour tout multi-indice $|L|=n-p$ on désigne par $R:=\complement L$ le multi-indice complémentaire de $L$ dans l'ensemble $\{1,..,n\}$. On a alors que le courant
$$
\Theta _{L,L}\cdot \omega _n= i^{p^2}2^{-n} \Theta\wedge\zeta _R^*\wedge\bar{\zeta}_R^*
$$
peut être identifié avec une mesure de Radon positive sur l'ouvert $U$. Nous reprenons maintenant un calcul fait par Demailly dans \cite{Dem-1}, chapitre III.
On désigne par $R:=\complement K,\;Q:=\complement H$ les multi-indices complémentaires de $K$ et $H$ dans l'ensemble $\{1,..,n\}$ et avec $\varepsilon _{\bullet}:=\pm 1,\;\pm i$. Avec ces notations on aura alors:
\begin{eqnarray}\label{Teta-K,H}
&\displaystyle{ 
\Theta _{K,H}\cdot \omega _n=\pm i^{p^2}2^{-n} \Theta\wedge\zeta _R^*\wedge\bar{\zeta}_Q^*=}&\nonumber
\\\nonumber
\\
&\displaystyle{
=2^{-n}\Theta \wedge \sum_{a\in(\Z/4\Z)^p}\,\varepsilon _a\bigwedge_{1\leq s\leq p}\frac{i}{4}  (\zeta ^*_{r_s}+i^{a_s}\zeta^*_{q_s})\wedge (\overline{\zeta ^*_{r_s}+i^{a_s}\zeta^*_{q_s}}).}&
\end{eqnarray}
En effet il suffit de remarquer l'identité extérieure
\begin{eqnarray}\label{fund-posit}
4\zeta^*_j\wedge \bar{\zeta}_k^*=(\zeta^*_j+\zeta^*_k)\wedge\overline{(\zeta^*_j+\zeta^*_k)}-(\zeta^*_j-\zeta^*_k)\wedge\overline{(\zeta^*_j-\zeta^*_k)}\nonumber
\\\nonumber
\\
+i(\zeta^*_j+i\zeta^*_k)\wedge\overline{(\zeta^*_j+i\zeta^*_k)}-i(\zeta^*_j-i\zeta^*_k)\wedge\overline{(\zeta^*_j-i\zeta^*_k)}.
\end{eqnarray}
Le fait que les formes
$$
\bigwedge_{1\leq s\leq p}\frac{i}{4}  (\zeta ^*_{r_s}+i^{a_s}\zeta^*_{q_s})\wedge (\overline{\zeta ^*_{r_s}+i^{a_s}\zeta^*_{q_s}})
$$ 
sont fortement positives entraîne, par hypothèse, que les courants $\Theta _{K,H}\cdot \omega _n$ peuvent être identifié avec une mesure de Radon complexe sur l'ouvert $U$, ce qui montre que le courant $\Theta $ est d'ordre zéro. D'autre part on a les égalités 
\begin{eqnarray*}
&\displaystyle{
2^{-n}\Theta \wedge\bigwedge_{1\leq s\leq p}\Big(\sum_{a_s\in(\Z/4\Z)}\, \frac{i}{4}  (\zeta ^*_{r_s}+i^{a_s}\zeta^*_{q_s})\wedge (\overline{\zeta ^*_{r_s}+i^{a_s}\zeta^*_{q_s}}) \Big)=}&
\\
\\
\\
&\displaystyle{
=2^{-n}\Theta \wedge \bigwedge_{1\leq s\leq p}  \,(i\zeta _{r_s}^*\wedge\bar{\zeta}_{r_s} ^*+i\zeta _{q_s}^*\wedge\bar{\zeta}_{q_s} ^*  )=}&
\\
\\
\\
&\displaystyle{
=2^{-n}\Theta \wedge \sum_{t\in E} \, i^{p^2}\,\zeta _{M_t} ^*\wedge\bar{\zeta}_{M_t} ^* = \sum_{t\in E} \,\Theta _{H_t,H_t}\cdot \omega _n}&
\end{eqnarray*}
où $E$ est un ensemble d'indices de cardinalité inférieure ou égale à $2^p$, $M_t\subset R\cup Q$ est un $p$-multi-indice et $H_t:=\complement M_t$. En utilisant l'expression \eqref{Teta-K,H} on obtient alors l'inégalité suivante:
\begin{eqnarray}\label{psh1} 
\sup_{\displaystyle
\scriptstyle   f\in {\cal D}^0(V,\C)  
\atop
\scriptstyle |f |\leq 1 } \left|  \,\int\limits_{V}\Theta _{K,H}\cdot f \omega _n \right| \leq 2^p\sum_{L\supset K\cap H}\, \int\limits_{V}\Theta _{L,L}\cdot \omega _n< +\infty
\end{eqnarray}
pour tout $V\subset U$ relativement compact dans $U$.
On remarque de plus que le courant $\Theta \wedge \omega _p$ s'écrit sous la forme
$$
\Theta \wedge \omega _p=2^{n-p} \Big(\sum_{|L|=n-p}\Theta _{L,L}  \Big)\cdot \omega _n.
$$
Le courant $\Theta \wedge \omega _p$ détermine  une mesure de Radon Positive $\|\Theta \|_{\omega}$ donnée explicitement par la formule
$$
\|\Theta \|_{\omega}(A):=\inf_{U \supset A} 
 \;\int\limits_{U}\Theta \wedge \omega _p
$$
pour tout sous-ensemble $A\subset X$ relativement compact. L'inégalité $\eqref{psh1} $ montre alors que les mesures de Radon complexes déterminées par les courants 
$\Theta _{K,H}\cdot \omega _n $ sont absolument continues par rapport à la mesure $\|\Theta \|_{\omega}$ restreinte à l'ouvert trivialisant $U$, ce qui prouve que la mesure de Radon $\|\Theta \|_{\omega}$ est une masse du courant $\Theta$. 
\\
Nous montrons maintenant que la forme $\theta _{\omega }(x)\in \Lambda ^{n-p,n-p}_{_J}\,T_{X,x} ^*$ est positive pour $\|\Theta \|_{\omega}$-presque tout $x\in X$. On désigne par 
$$
FP_p(\zeta)\subset {\cal E}(\Lambda ^{p,p}_{_J}\,T_X ^*)(U)
$$ 
l'ensemble des $(p,p)$-formes fortement positives à coeficients constants par rapport au repère $(\zeta _1,...,\zeta _n)$ et on considère  un sous-ensemble $(\varphi _{\nu} )_{\nu \in\N}\subset FP_p(\zeta) $ dense dans $FP_p(\zeta)$. Soit $\xi _{\omega}\in {\cal E}(\Lambda ^{n,n}_{_J}\,T_X)(X)$ le $(n,n)$-champ de vecteurs tel que $\omega _n(\xi _{\omega})=1$ sur $X$. On remarque que pour tout $(p,p)$-forme $\varphi  \in {\cal E}(\Lambda ^{p,p}_{_J}T_X^*)(X)$ et tout Borelien $A\subset X$ on a les identités
$$
\int\limits_{A}\Theta \wedge \varphi = \int\limits_{A}(\Theta \wedge \varphi)(\xi_{\omega})\cdot \omega _n= \int\limits_{A}(\theta _{\omega}\wedge \varphi )(\xi_{\omega})\,\|\Theta \|_{\omega}.
$$
On a alors que l'ensemble 
$$
E_{\nu }:=\{x\in \mbox{Dom}\,\theta _\omega\cap U \,|\,\theta _\omega (x)\wedge \varphi_{\nu } (x)<0 \}
$$
est un ensemble de $\|\Theta \|_{\omega}$-mesure nulle (ici $\mbox{Dom}\,\theta _\omega$ désigne le domaine du représentant $\theta _\omega$). Le fait que 
$$
\Lambda ^{p,p}_{_J}\,T_{X,x}^*=\{\varphi (x)\,|\,\varphi \in FP_p(\zeta)\}
$$
pour tout $x\in U$ combiné avec le fait que, par densité, pour tout  $\varphi \in FP_p(\zeta)$
 il existe une suite $(\nu _l)_l$ telle que $\varphi =\lim_{l\rightarrow +\infty}\varphi _{\nu _l}$ entraînent
$$
\theta _\omega (x)\wedge \varphi (x)=\lim_{l\rightarrow +\infty}\theta _\omega(x)\wedge \varphi _{\nu _l}(x) \geq 0 
$$ 
pour tout $x\in \mbox{Dom}\,\theta _\omega\cap U \setminus \cup_{\nu }E_{\nu}$. Ceci entraîne la conclusion voulue sur la forme $\theta _\omega$. 
\hfill $\Box$
\\
\\
Voyons maintenant quelques exemples fondamentaux de $(1,1)$-courant positif sur les variétés presque complexes.
\subsection{Exemples fondamentaux de courants positifs sur les variétés presque complexes.} 
On commençe par une définition.
\begin{defi} 
Une sous-variété $Y\subset X$ de dimension $2p$ d'une variété presque complexe $(X,J)$ est dite presque complexe si $J(T_Y)=T_Y$.
\end{defi}
Un exemple de sous-variété presque complexe est constitué par les images $\gamma(\proj^1_{_{\C} } )\subset X$ des courbes $J$-holomorphes régulières. Les résultats qui suivront vont assurer l'existence d'exemples de sous-variétés presque complexes de dimension complexe supérieure à un. On commence par rappeler la proposition suivante (voir \cite{McD-Sa}).
\begin{prop} Soit $X$ une variété différentielle de dimension réelle $2n$. S'il existe une $2$-forme $\omega \in {\cal E}(\Lambda ^2T_X^*)(X)$ non dégénéré alors l'espace des structures presque complexes compatibles avec  $\omega$
$$
{\cal J}_{X,\omega }:=\Bigl\{J\in {\cal E}(T^*_X\otimes_{_{\R}} T_X)(X)\,|\,J^2=-\I,\,\omega (Ju,Jv)=\omega (u,v),\,\omega (u,Ju)>0\,\forall u,v\in T_{X,x}\smallsetminus \left\{0_x\right\}  \Bigr\}   
$$
est non vide et contractile.
\end{prop} 
On a aussi la proposition suivante, (voir l'article de Audin dans l'ouvrage \cite{Ad}).
\begin{prop}
Soit $X$ une variété différentielle de dimension réelle $2n$ admettant une $2$-forme $\omega \in {\cal E}(\Lambda ^2T_X^*)(X)$ non dégénérée et soit $Y\subset X$ une sous-variété telle que $i_{Y}^*\omega$ soit non dégénérée. Il existe alors une structure presque complexe $J\in {\cal J}_{X,\omega }$ telle que $(Y,J_{|_Y})$ soit une sous-variété presque complexe de $(X,J)$.
\end{prop} 
On a le résultat fondamental suivant du à S.K Donaldson (voir \cite{Don}).
\begin{theo}\label{Donaldson} 
Soit $(X,\omega)$ une variété symplectique compacte de dimension réelle $2n$. Pour tout $p=1,...,n$ il existe des sous-variétés symplectiques $(Y_p,i_{Y_p}^*\omega)$ fermées de dimension réelle $2p$.
\end{theo} 
Le premier exemple de courant positif qu'on considère est le courant d'intégration sur une sous-variété presque complexe $Y$ de dimension $2p$ de mesure localement finie avec l'orientation canonique donnée par la structure presque complexe $J_{|_Y}\in {\cal E}(T^*_Y\otimes_{_{\R}} T_Y)$. Le courant $[Y]\in {\cal D}'_{2p}(X)$ s'identifie naturellement avec un élément de l'espace  ${\cal D}'_{p,p}(X)$  étant $\int_Y \varphi =\int_Y \varphi^{p,p}$ pour tout $\varphi \in {\cal D}^{2p}(X)$, où 
$\varphi^{p,p}$ désigne la composante de type $(p,p)$ de la forme $\varphi$. Le courant $[Y]$ est évidemment positif grâce à la  propriété $\ref{corrpos}.c$. Sous les hypothèses du théorème $\ref{Donaldson}$ on a alors l'existence de courants $[Y_p]$ lesquels sont de bidegré $(n-p,n-p)$ et positifs par rapport à une structure presque complexe $J_p\in {\cal J}_{X,\omega }$. De plus les courants en question sont fermés, i.e $d[Y_p]=0$ en conséquence de la formule de Stokes. La notion intuitive de la masse $\|[Y]\|_{\omega}$ est clarifiée par le lemme suivant qui est une généralisation immédiate d'un résultat bien connu dans le cas des variétés complexes (voir le chapitre III dans l'ouvrage de Demailly \cite{Dem-1}).
\begin{iw} Soit $Y\subset X$ une sous-variété orientable et orientée de dimension $2p$ d'une variété presque complexe $(X,J)$ munie d'une métrique $\omega \in {\cal E}(\Lambda ^{1,1}_{_J}T_X^*)(X)$. Si on désigne par $dV_{Y,\omega }$ la forme de volume associée à la restriction à $T_Y$ de la métrique riemannienne 
$g(\cdot,\cdot):=\omega(\cdot,J\cdot) $ associée à la métrique $\omega $ on a l'existence d'une fonction $\alpha \in {\cal C}^0(Y,[-1,1]) $ telle que 
$\omega _{p_{|_Y}}=\alpha \cdot dV_{Y,\omega }$. De plus $|\alpha |=1$ si et seulement si $Y$ est une sous-variété presque complexe. Dans ce cas $\alpha =1$ si l'orientation de $Y$ coïncide avec l'orientation canonique donnée par la structure $J_{|_Y} $ et $\alpha =-1$ sinon. La fonction $\alpha $ est identiquement nulle si et seulement si $Y$ est une sous-variété $\omega $-isotropique.
\end{iw} 
Le théorème suivant nous fournit un autre exemple fondamental de $(1,1)$-courant positif.
\begin{theo}\label{pluri-positif}
Soit $(X,J)$ une variété presque complexe connexe et $f\in Psh(X,J)$. Alors ou bien $f\equiv -\infty$ ou bien $f\in L^1_{loc}(X)$. Dans ce dernier cas le $(1,1)$-courant $i\partial_{_J }\bar{\partial}_{_J }f$ est positif.
\end{theo}   
$Preuve$
\\
{\bf Intégrabilité locale de $f$}. Avec les notations du théorème \ref{Jplat} on a que pour tout $x\in U_{x_0}$ 
l'application 
\begin{diagram}[height=1cm,width=1cm]
&\varphi _x:(0,\delta )\times S^1\times  S^{2n-1}(T_{X,x} ) &\rTo&
X
\\
&(r,\theta,v)&\longmapsto& \Phi(re^{i\theta},v)& 
\end{diagram}
est une submersion de classe $\ci$. 
Par hypothèse on a l'inégalité de la moyenne
$$
f(x)\leq \frac{1}{2\pi r}  \int\limits_0^{2\pi}f\circ  \varphi _x(r,\theta ,v)\,d\theta. 
$$ 
On considère les ouverts relativement compacts
$$
C_{r_1,r_2} (x):=\varphi _x\Big((r_1,r_2 )\times S^1\times  S^{2n-1}(T_{X,x} )\Big), \;0<r_1<r_2<\delta 
$$
et la forme de volume de classe $\ci$ 
$$
d{\cal V}_x(p):=\int\limits_{(r,\theta ,v)\in \varphi _x^{-1}(p) } dr\,d\theta \,d\sigma (v)$$
sur l'ouvert $\mbox{Im}\varphi _x$, où $d\sigma$ désigne la forme de volume sur la sphère $S^{2n-1}(T_{X,x} )$.
Avec ces notations on a alors
\begin{eqnarray}\label{ineqSH}
\int\limits_{p\in C_{r_1,r_2} (x)}f(p)\,d{\cal V}_x(p)=
\int\limits_{v\in S^{2n-1}(T_{X,x})}d\sigma (v)\int\limits_{r_1}^{r_2}dr \int\limits_0^{2\pi}f\circ  \varphi _x(r,\theta ,v)\,d\theta\geq f(x)\,K_{r_1,r_2}     
\end{eqnarray}  
où $K_{r_1,r_2}>0$ est une constante.
Soit $W\subset X$ l'ensemble des points $p\in X$ tels que la fonction $f$ soit intégrable sur un voisinage de $p$. Par définition le sous-ensemble $W$ est ouvert en $X$ et $f>-\infty$ presque partout sur $W$. Si $p\in \overline{W}$, on peut choisir un point $x\in W$ tel que 
$f(x)>-\infty$ et $p\in C_{r_1,r_2} (x)$. On déduit alors d'après l'inégalité \eqref{ineqSH}, que la fonction $f$ est intégrable sur le voisinage $C_{r_1,r_2} (x)$ de $p$, ce qui montre que $p\in W$ et donc que $W$ est aussi fermé en $X$. On a alors soit $W=X$, soit $W=\emptyset$. Dans le dernier cas l'inégalité \eqref{ineqSH}, implique $f\equiv -\infty$. On a donc prouvé que soit $f\equiv -\infty$ soit $f\in L^1_{loc}(X)$. 
\\
\\
{\bf Positivité du courant $i\partial_{_J }\bar{\partial}_{_J }f$}. On montre d'abord que pour tout 
$\xi\in P_{_J}(U_{x_0} ,T_X)$ la distribution $i\partial_{_J }\bar{\partial}_{_J }f\,(\xi ,J\xi )$ est positive sur $U_{x_0}$. Pour tout $x\in U_{x_0}$ soient 
$(U_x,\sigma _{\xi }^{-1})$
$$
\sigma_{\xi}  :B^1_{\delta}\times B^{n-1}_{\delta  }
\longrightarrow U_x\subset U_{x_0}
$$
les coordonnées du lemme $\ref{locJ-plat}$. En rappelant l'expression explicite $\eqref{explicit(1,1)cour}$ du courant  $i\partial_{_J }\bar{\partial}_{_J }f$ on aura pour tout $\xi \in {\cal E}(T_X)(U_{x_0} )$ les égalités suivantes:
\begin{eqnarray*} 
&i\partial_{_J }\bar{\partial}_{_J }f\,(\xi ,J\xi )=2\partial_{_J }\bar{\partial}_{_J }f\,(\xi^{1,0}  ,\xi^{0,1})
 =2(\xi ^{1,0}.\,\xi ^{0,1}.\,f-[\xi^{1,0}  ,\xi^{0,1}]^{0,1} .\,f)=&
\\
\\
&= \displaystyle{ \frac{1}{2}}(\xi .\,\xi .\,f+ J\xi .\,J\xi .\,f+J[\xi ,J\xi ].\,f )   &
\end{eqnarray*}
Le fait que dans notre cas $[\xi ,J\xi ]=0$, implique les expressions:
\begin{eqnarray*} 
i\partial_{_J }\bar{\partial}_{_J }f\,(\xi ,J\xi )=2\xi ^{1,0}.\,\xi ^{0,1}.\,f=\frac{1}{2}(\xi .\,\xi .\,f+ J\xi .\,J\xi .\,f )=
\frac{1}{2}(\sigma _{\xi }^{-1})^*\Delta_{z_1}(f\circ \sigma _{\xi })   
\end{eqnarray*}
où $\Delta_{z_1}:=\partial^2_{x_1}+\partial^2_{y_1} $ désigne le Laplacien par rapport à la variable $z_1=x_1+i\,y_1\in B^1_{\delta}$ dans l'ouvert $B^1_{\delta}\times B^{n-1}_{\delta  }$. Grâce au théorème de Fubini on en déduit l'inégalité
\begin{eqnarray*}
\int\limits_{B^1_{\delta}\times B^{n-1}_{\delta  } }(f\circ \sigma _{\xi })\cdot \Delta_{z_1}\varphi \,d\lambda \geq 0
\end{eqnarray*}
pour tout $\varphi \in{\cal D}(B^1_{\delta}\times B^{n-1}_{\delta  }),\;\varphi \geq 0$. Le Laplacien 
$\Delta_{z_1}(f\circ \sigma _{\xi }) $ est donc positif, ce qui prouve la positivité de la distribution $i\partial_{_J }\bar{\partial}_{_J }f\,(\xi ,J\xi )$ sur l'ouvert $U_{x_0}$ pour tout champ $\xi\in P_{_J}(U_{x_0} ,T_X)$. Nous montrons maintenant que le courant $i\partial_{_J }\bar{\partial}_{_J }f$  est d'ordre zéro. 
\\
Soit $\zeta_1,...,\zeta_n$ un repère complexe du fibré des $(1,0)$ vecteurs tangents  $T^{1,0}_{X,J|U_{x_0} }$. On déduit d'après l'identité extérieure \eqref{fund-posit}, (avec $\zeta$ à la place de $\zeta ^*$) l'existence de champs de vecteurs $\rho _k\in {\cal E}(T^{1,0}_{X,J})(U_{x_0}),\,k=1,...,n^2$ du type 
$\rho _k=\zeta _{s_k}+i^{a_l}\zeta _{t_k}, \,a\in \Z/4\Z$ tels que les $(1,1)$-champs de vecteurs $(\rho _k\wedge \bar{\rho}_k)_{k=1}^{n^2}$ forment un repère complexe du fibré $\Lambda ^{1,1}_{_J}\,T_{U_{x_0} } $. On choisit un point $x\in U_{x_0}$ et $\xi_k\in P_{_J}(U_{x_0}  ,T_X)$  tels que $\xi^{1,0} _k (x)=\rho_k(x)$. On aura alors que les $(1,1)$-champs de vecteurs $(\xi ^{1,0} _k\wedge \xi ^{0,1} _k)_{k=1}^{n^2}$ forment un repère complexe du fibré $\Lambda ^{1,1}_{_J}\,T_{V_x}$ ou $V_x\subset U_{x_0}$ est un voisinage ouvert du point $x$. L'identité 
$$
i\partial_{_J }\bar{\partial}_{_J }f\,(\xi_k ,J\xi_k )=2\partial_{_J }\bar{\partial}_{_J }f\,(\xi ^{1,0} _k\wedge \xi ^{0,1} _k)
$$ 
montre que le courant $i\partial_{_J }\bar{\partial}_{_J }f$  est d'ordre zéro. 
\\
Venons-en maintenant à la positivité du courant en question. Soit $\mu$ une masse du courant $i\partial_{_J }\bar{\partial}_{_J }f$ et considerons l'écriture $i\partial_{_J }\bar{\partial}_{_J }f=\theta\cdot \mu$. Nous montrons que la forme $\theta (x)\in \Lambda ^{1,1}_{_J}\,T_{X,x} ^*$ est positive pour $\mu$-presque tout $x\in X$. On désigne par 
$$
\Q T_{X|U_{x_0}}:=T_{X|U_{x_0}}\cap (\Q^{2n}\times\Q^{2n}),
$$
en supposant que l'ouvert 
$U_{x_0}$ est un ouvert coordonné. D'après la preuve du théorème $\ref{Jplat} $ il existe une famille dénombrable de champs $(\xi _{\nu})_{\nu \in \N}\subset P_{_J}(U_{x_0}  ,T_X)$ telle que pour tout $v\in\Q T_{X|U_{x_0}}\smallsetminus 0_X$ il existe $\nu \in \N$ tel que $\xi _{\nu}(\pi (v))=v$ ($\pi$ désigne la projection sur le fibré tangent) et pour tout $x\in U_{x_0}$ l'ensemble $(\xi _{\nu}(x))_{\nu \in \N}$ est dense dans $T_{X,x}$. La positivité de la distribution $i\partial_{_J }\bar{\partial}_{_J }f\,(\xi ,J\xi )$ sur l'ouvert $U_{x_0}$entraîne que l'ensemble 
$$
E_{\nu }:=\{x\in \mbox{Dom}\,\theta \cap U_{x_0}\,|\,\theta (\xi _{\nu},J\xi _{\nu})(x)<0 \}
$$
est de $\mu$-mesure nulle (ici $\mbox{Dom}\,\theta$ désigne le domaine du représentant $\theta$). Pour tout $x\in \mbox{Dom}\,\theta \cap U_{x_0} \setminus \cup_{\nu }E_{\nu}$ et pour tout $v\in T_{X,x}$ considerons une suite $(\nu _l)_l$ telle que $v=\lim_{l\rightarrow +\infty}\xi _{\nu_l}(x)$. La limite
$$
\theta (v,Jv)(x)=\lim_{l\rightarrow +\infty}\theta(\xi _{\nu_l},J\xi _{\nu_l})(x )\geq 0 
$$ 
entraîne alors la conclusion voulue sur la forme $\theta$. \hfill $\Box$
\section{Les potentiels des courants positifs de type $(1,1)$ sur les  
\\
variétés presque complexes }
Dans cette section nous proposons une conjecture réciproque du théorème $\ref{pluri-positif}$ qu'on énonce sous la forme suivante.
\begin{mut}\label{cas-Gen-Psh} 
Soit $(X,J)$ une variété presque complexe de dimension complexe $n$ et soit $u\in {\cal D}'_{2n}(\R)(X)$ une distribution réelle telle que le $(1,1)$-courant 
$i\partial_{_J }\bar{\partial}_{_J }u\in{{\cal D}'}^{1,1} (X)$ soit positif. Alors il existe une unique fonction $f\in Psh(X,J)\cap L^1_{loc}(X)$ telle que la distribution correspondante coïncide avec la distribution $u$.
\end{mut}
Il est bien connu que la conjecture est vraie dans le cas complexe intégrable (voir \cite{Dem-1}). 
\\
{\bf Remarque 1}. On considère l'opérateur 
$$
d^c_{_{J}}:=\frac{i}{2}(\bar{\partial}_{_J }-\partial_{_J }).
$$
En degré zéro il se réduit à la forme $d^c_{_{J}}:=-\frac{1}{2}df\circ J$. En utilisant les identités fondamentales de la géométrie presque complexe on déduit facilement qu'en degré zéro on a l'identité 
$$
i\partial_{_J }\bar{\partial}_{_J }=  dd^c_{_{J}}+i\theta _{_{J} }\bar{\partial}_{_J }-i\bar{\theta}_{_{J}}\partial_{_J }, 
$$
qui montre de quelle façon la torsion de la structure presque complexe représente l'obstruction pour le $(1,1)$-courant $i\partial_{_J }\bar{\partial}_{_J }u$ à être $d$-fermé.
D'après l'identité précédente on déduit alors l'égalité 
$$
dd^c_{_{J}}\,u\,(\xi ,J\xi )=i\partial_{_J }\bar{\partial}_{_J } u\,(\xi).
$$
On a donc que le
$(1,1)$-courant $i\partial_{_J }\bar{\partial}_{_J }u$ est positif si et seulement si pour tout champs de vecteurs réel $\xi $ la distribution $dd^c_{_{J}}\,u(\xi ,J\xi )$ est positive. On remarque de plus que comme dans le cas complexe intégrable, (cf. \cite{Dem-1}) on a d'après la formule de Stokes l'égalité
\begin{eqnarray}
\int\limits_U\varphi \wedge dd^c_{_{J}}\psi-dd^c_{_{J}}\varphi \wedge \psi=
\int\limits_{\partial U}\varphi \wedge d^c_{_{J}}\psi-d^c_{_{J}}\varphi \wedge \psi
\end{eqnarray}  
pour tout ouvert $U\subset X$ relativement compact  à bord ${\cal C}^1$ par morceaux et pour tout $\varphi \in {\cal C}^2(\Lambda ^{p,p}_{_J}T_X^*)(\overline{U})$ et $\psi\in {\cal C}^2(\Lambda ^{q,q}_{_J}T_X^*)(\overline{U}), \;p+q=n-1$. En utilisant la formule précédente et le fait que $\partial_{_J }\bar{\partial}_{_J }+\bar{\partial}_{_J }\partial_{_J }=0$ en bidegré $(n-1,n-1)$, on déduit pour tout $\varphi \in {\cal D}^{n-1,n-1}(X)$ les égalités suivantes
$$
\int\limits_X i\partial_{_J }\bar{\partial}_{_J }u\wedge \varphi =
\int\limits_X  dd^c_{_{J}}u\wedge \varphi =
\int\limits_X u \cdot dd^c_{_{J}}\varphi =
\int\limits_X u \cdot i\partial_{_J }\bar{\partial}_{_J }\varphi.
$$ 
\\
{\bf Remarque 2}. Soit $\omega \in {\cal E}(\Lambda ^{1,1}_{_J}T_X^*)(X)$ une métrique hermitienne sur $T_{X,J}$. On définit le Laplacien de $u\in {\cal D}'_{2n}(\R)(X)$ par rapport à la structure presque complexe $J$ et la métrique $\omega$ par la formule
$$
\Delta _{_{J,\omega}}u:=\mbox{Trace}_{_{\omega} } (i\partial_{_J }\bar{\partial}_{_J }u)=
\frac{ n\cdot i\partial_{_J }\bar{\partial}_{_J }u\wedge \omega ^{n-1} }{2\cdot\omega ^n}.     
$$
Soit $(\xi _k)_k\in {\cal E}(T_{X,J})(U)^{\oplus n}$ un repère local complexe $\omega $-orthonormé du fibré $T_{X,J}$ et soit $\zeta _k:=\xi^{1,0}_k$. En rappelant l'écriture locale  $\eqref{exploccour1,1}$ on obtient l'égalité
\begin{eqnarray*}
&\displaystyle{
\Delta _{_{J,\omega}}u =\sum_{k=1}^n (\zeta _k\,.\bar{\zeta}_k .\,u-[\zeta _k,\bar{\zeta}_k]^{0,1}.\,u )=}&
\\
\\
&\displaystyle{
=\frac{1}{4} \sum_{k=1}^n (\xi _k\,.\,\xi _k\,.u+J\xi _k\,.\,J\xi _k\,.u+J[\xi _k,J\xi _k]\,.u)= \frac{1}{2} \sum_{k=1}^ni\partial_{_J }\bar{\partial}_{_J }u(\xi _k,J\xi _k)}.&
\end{eqnarray*}
Ce calcul montre que le symbole de $\Delta _{_{J,\omega}}$ coïncide avec le symbole du Laplacien classique $\Delta\equiv\Delta _{_{J_0,\omega_0}}$ sur $\C^n$, où
$\omega_0=\frac{i}{2}\partial_{_{J_0} }\bar{\partial}_{_{J_0}}|z|^2$ est la métrique $J_
0$-invariante plate sur $\C^n$. On obtient alors que l'opérateur de Green de $\Delta _{_{J,\omega}}$ coïncide avec l'opérateur de Green classique de $\Delta$ au sens des opérateurs pseudodifférentiels.
Si le courant $i\partial_{_J }\bar{\partial}_{_J }u$ est positif alors $\Delta _{_{J,\omega}}u$ est une mesure de Radon positive. On déduit alors d'après la théorie classique des opérateurs elliptiques d'ordre deux que $u\in W^{1,1}_{loc}(X):= \{v\in L^1_{loc}(X)\,|\,dv\in L^1_{loc}(T^*_X)(X)\} $, (cf. \cite{Sta}, paragraphe 9, théorèmes 9.1 et 9.4).
\\
Nous montrons maintenant la conjecture dans le cas particulier suivant.
\begin{theo}\label{caspart} 
Soit $(X,J)$ une variété presque complexe et $f:X\longrightarrow [-\infty,+\infty)$ une fonction semi-continue supérieurement telle que $f$ soit continue sur l'ensemble $X\smallsetminus f^{-1}(-\infty)$, $f\in L^1_{loc}(X)$ et telle que le $(1,1)$-courant 
$i\partial_{_J }\bar{\partial}_{_J }f\in{{\cal D}'}^{1,1} (X)$ soit positif. Alors $f\in Psh(X,J)$.
\end{theo}   
Avant de passer à la preuve du théorème $\ref{caspart}$ nous aurons besoin de quelques notions et résultats préliminaires.
On commence par prouver le lemme suivant.
\begin{lem}\label{red-cont}  
Sous les hypothèses du théorème $\ref{caspart}$, le courant $i\partial_{_J }\bar{\partial}_{_J }\log (e^f+\varepsilon )$ est positif pour tout $\varepsilon >0$.
\end{lem} 
$Preuve$. Le fait que $f\in L^1_{loc}(X)$ implique que l'intérieur de l'ensemble $f^{-1}(-\infty)$ est vide. On déduit que pour tout $x\in f^{-1}(-\infty)$ et pour tout $r>0$ il existe $y\in B^n_r(x)\smallsetminus f^{-1}(-\infty)$.
L'hypothèse de semi-continuité combinée avec la continuité de la fonction $f$ sur l'ensemble $X\smallsetminus f^{-1}(-\infty)$ entraînent alors la continuité de la fonction $e^f$ sur tout $X$, (en particulier l'ensemble $f^{-1}(-\infty)$ est fermé dans $X$). On obtient alors la continuité des fonctions $f_{\varepsilon} :=\log (e^f+\varepsilon )$. On remarque que si $u$ est une fonction de classe ${\cal C}^2$ on a la formule
$$
i\partial_{_J }\bar{\partial}_{_J }u_{\varepsilon}=\frac{e^u}{e^u+\varepsilon} \, i\partial_{_J }\bar{\partial}_{_J }u+
\frac{\varepsilon e^u}{(e^u+\varepsilon)^2}\,i\partial_{_J }u\wedge\bar{\partial}_{_J }u 
$$
et la $(1,1)$-forme $i\partial_{_J }u\wedge\bar{\partial}_{_J }u$ est positive. En effet pour tout champs de vecteurs réels $\xi$ on a les égalités
\begin{eqnarray*}
&\displaystyle{i\partial_{_J }u\wedge\bar{\partial}_{_J }u\,(\xi ,J\xi )=i\partial_{_J }u(\xi )\cdot\bar{\partial}_{_J }u(J\xi )-i\partial_{_J }u(J\xi )\cdot\bar{\partial}_{_J }u(\xi )=}&
\\
\\
&\displaystyle{=2\partial_{_J }u(\xi )\cdot\bar{\partial}_{_J }u(\xi )=\frac{1}{2}(du(\xi )^2+du(J\xi )^2)\geq 0}&
\end{eqnarray*}  
On en déduit alors que si notre fonction $f$ est de classe ${\cal C}^2$ la $(1,1)$-forme $i\partial_{_J }\bar{\partial}_{_J }f_{\varepsilon}$ est positive. Dans le cas général nous considérons une famille de noyaux régularisants $(\rho _{\eta} )_{\eta>0} $ sur un ouvert coordonnée $V\subset X$ et les fonctions $f^{\eta}:=f*\rho _{\eta},\,f^{\eta} _{\varepsilon} :=\log (e^{f^{\eta} } +\varepsilon )\in{\cal E}(U,\R)$ où $U\subset V$ est un ouvert relativement compact dans $V$. (On remarque que si la structure presque complexe n'est pas intégrable les $(1,1)$-formes $i\partial_{_J }\bar{\partial}_{_J }f^{\eta}$ ne sont pas positives). Pour prouver la positivité du courant $i\partial_{_J }\bar{\partial}_{_J }f_{\varepsilon}$ on remarque que pour tout forme $\varphi \in{\cal D}^{2n}(U)$ positive et pour tout champ de vecteurs réels $\xi \in {\cal E}(T_X)(U)$ on a les égalités suivantes
\begin{eqnarray*}
&\displaystyle{
\int\limits_U i\partial_{_J }\bar{\partial}_{_J }f_{\varepsilon}(\xi ,J\xi )\,\varphi =
\lim_{\eta \rightarrow 0} \int\limits_U i\partial_{_J }\bar{\partial}_{_J }f^{\eta} _{\varepsilon}(\xi ,J\xi )\,\varphi=}&
\\
\\
&\displaystyle{=
\lim_{\eta \rightarrow 0} \left[\int\limits_U \frac{e^{f^{\eta} }}{e^{f^{\eta} }+\varepsilon} \, i\partial_{_J }\bar{\partial}_{_J }f^{\eta}(\xi ,J\xi )\,\varphi+
\int\limits_U\frac{\varepsilon e^{f^{\eta} }}{(e^{f^{\eta} }+\varepsilon)^2}\,i\partial_{_J }f^{\eta}\wedge\bar{\partial}_{_J }f^{\eta}(\xi ,J\xi )\,\varphi\right]}.& 
\end{eqnarray*} 
 De plus on va montrer l'égalité
\begin{eqnarray}\label{machedevofa}  
\lim_{\eta \rightarrow 0} \int\limits_U \frac{e^{f^{\eta} }}{e^{f^{\eta} }+\varepsilon} \, i\partial_{_J }\bar{\partial}_{_J }f^{\eta}(\xi ,J\xi )\,\varphi=
\int\limits_U\frac{e^f}{e^f+\varepsilon} \, i\partial_{_J }\bar{\partial}_{_J }f(\xi ,J\xi )\,\varphi\geq 0.
\end{eqnarray} 
On aura alors
$$
\int\limits_U i\partial_{_J }\bar{\partial}_{_J }f_{\varepsilon}(\xi ,J\xi )\,\varphi =
\int\limits_U\frac{e^f}{e^f+\varepsilon} \, i\partial_{_J }\bar{\partial}_{_J }f(\xi ,J\xi )\,\varphi+\lim_{\eta \rightarrow 0} \int\limits_U\frac{\varepsilon e^{f^{\eta} }}{(e^{f^{\eta} }+\varepsilon)^2}\,i\partial_{_J }f^{\eta}\wedge\bar{\partial}_{_J }f^{\eta}(\xi ,J\xi )\,\varphi.
$$
La dernière limite est positive car la forme $i\partial_{_J }f^{\eta}\wedge\bar{\partial}_{_J }f^{\eta}$ est positive.
Pour prouver l'égalité $\eqref{machedevofa}$ il suffit de montrer que la suite $(i\partial_{_J }\bar{\partial}_{_J }f^{\eta})_{\eta>0}$ converge vers le courant $i\partial_{_J }\bar{\partial}_{_J }f$ aussi dans la topologie faible des courants d'ordre zéro. Ce fait combiné avec le fait que la suite de fonctions 
${e^{f^{\eta} }}/(e^{f^{\eta} }+\varepsilon)$ converge uniformément vers la fonction $e^f/(e^f+\varepsilon )$ prouve l'égalité $\eqref{machedevofa}$. Pour prouver la convergence de la suite $(i\partial_{_J }\bar{\partial}_{_J }f^{\eta})_{\eta>0}$ dans la topologie faible des courants d'ordre zéro il suffit de montrer la convergence de la suite 
$$
(i\partial_{_J }\bar{\partial}_{_J }f^{\eta}-(i\partial_{_J }\bar{\partial}_{_J }f)*\rho _{\eta})_{\eta>0} 
$$ 
dans la même topologie étant donné que la suite $((i\partial_{_J }\bar{\partial}_{_J }f)*\rho _{\eta})_{\eta>0}$ est convergente dans cette topologie. D'après le lemme $\ref{critconvzero}$ il suffit donc de remarquer l'inégalité 
\begin{eqnarray*}
\sup_{\eta>0}\mu(i\partial_{_J }\bar{\partial}_{_J }f^{\eta}-(i\partial_{_J }\bar{\partial}_{_J }f)*\rho _{\eta})(U)\leq C\|f\|_{W^{1,1}(U)}<\infty
\end{eqnarray*}
qui découle du lemme de K.O. Friedrichs (cf. \cite{Hor-1}). \hfill $\Box$
\\
\\ 
On rappelle la définition suivante.
\begin{defi}
Un sous-ensemble $A\subset U$ d'un ouvert $U\subset \R^m$ est dit polaire si pour tout point $x\in U$ il existe un voisinage ouvert connexe $V_x\subset U$ de $x$ et une fonction $u$ sous-harmonique sur $V_x$, $u\not \equiv -\infty$, telle que $A\cap V_x\subset  \{y\in V_x\,|\, u(y)=-\infty\}$. 
\end{defi}  
D'après le théorème \ref{equiv-SH} on a qu'un sous-ensemble polaire est de mesure de Lebesgue nulle.  
On a le théorème classique suivant (cf. \cite{Dem-1}, chapitre I).
\begin{theo}\label{psh-exten}
Soit $A\subset U$ un sous-ensemble polaire fermé et soit $v$ une fonction sous-harmonique sur l'ouvert $U\smallsetminus A$, borné supérieurement sur un voisinage de tout point de $A$. Il existe alors une unique extension sous-harmonique $\tilde{v}$ de $v$ sur $U$. En particulier si $v$ est une fonction continue sur $U$ et sous-harmonique sur l'ouvert $U\smallsetminus A$ alors $v$ est sous-harmonique sur $U$.
\end{theo}  
{\bf Preuve du théorème \ref{caspart}}. D'après le lemme \ref{red-cont} 
il suffit de montrer le théorème dans le cas d'une fonction continue étant donné que la fonction $f$ est limite décroissante des fonctions continues $f_{\varepsilon} :=\log (e^f+\varepsilon )$ (lorsque $\varepsilon $ tend vers zéro) et une limite décroissante de fonctions plurisousharmoniques est plurisousharmonique. A partir de maintenant on suppose donc $f$ continue et on remarque que pour tout courbe $J$-holomorphe $\gamma :B^1_{\rho}\longrightarrow X$ la fonction $f\circ \gamma$ est sous-harmonique sur $B^1_{\rho}$ si et seulement si elle est sous-harmonique sur l'ouvert 
$$
\{z\in B^1_{\rho}\,|\, d_z\gamma\not=0\}. 
$$
Ceci découle du fait que l'ensemble $\{z\in B^1_{\rho}\,|\, d_z\gamma=0\}$ est fini (voir \cite{McD-1}, chapitre II) et du théorème \ref{psh-exten}. On peut donc supposer que la  courbe $J$-holomorphe $\gamma :B^1_{\rho}\longrightarrow X$ est un plongement.
\\
D'autre part on remarque qu'une fonction continue $u:\Omega \subset \R^n\longrightarrow \R$ est sous-harmonique si et seulement si $\Delta u\geq 0$. En effet en considérant une famille de noyaux régularisants usuels $(\rho _{\varepsilon } )_{\varepsilon >0}$ on a $0\leq (\Delta u)*\rho _{\varepsilon}=\Delta (u*\rho _{\varepsilon})$. On déduit alors l'inégalité
$$
\frac{1}{\lambda(B_r(x))}   \int\limits_{B_r(x) }u*\rho _{\varepsilon} \,d\lambda \geq (u*\rho _{\varepsilon})(x).
$$ 
En passant à la limite pour $\varepsilon$ tendent vers zéro on déduit la sous-harmonicité de $u$. On va donc montrer l'inégalité $\Delta (f\circ \gamma)\geq 0$ pour tout plongement $J$-holomorphe $\gamma :B^1_{\rho}\longrightarrow X$. D'après le lemme $\ref{perturbCurbJ-hol}$ on a, quitte à restreindre $\rho >0$, l'existence d'un plongement
$$
\sigma :B^1_{\rho }\times B^{n-1}_{\rho}\longrightarrow \sigma (B^1_{\rho }\times B^{n-1}_{\rho})\subset X
$$
qui préserve les orientations canoniques tel que $\sigma (\cdot, z_2)$ soit une courbe $J$-holomorphe pour tout $z_2\in B^{n-1}_{\rho}$ et $\sigma (z_1,0)=\gamma (z_1),\,z_1=t+is$. 
Le fait que le champ $\xi:=d\sigma (\frac{\partial}{\partial t})\circ \sigma ^{-1}$ soit $J$-plat sur l'ouvert $\sigma (B^1_{\rho }\times B^{n-1}_{\rho})$ implique les égalités
\begin{eqnarray*} 
i\partial_{_J }\bar{\partial}_{_J }f\,(\xi ,J\xi )=2\xi ^{1,0}.\,\xi ^{0,1}.\,f=\frac{1}{2}(\xi .\,\xi .\,f+ J\xi .\,J\xi .\,f )=
\frac{1}{2}(\sigma^{-1})^*\Delta_{z_1}(f\circ \sigma)   
\end{eqnarray*}
où $\Delta_{z_1}:=\partial^2_{t}+\partial^2_{s}$ désigne le Laplacien par rapport à la variable $z_1=t+i\,s\in B^1_{\rho}$ dans l'ouvert $B^1_{\rho }\times B^{n-1}_{\rho}$. 
Le fait que le plongement $\sigma $ préserve les orientations canoniques implique l'inégalité 
$$
\Delta_{z_1}(f\circ \sigma)\geq 0
$$ 
sur l'ouvert $B^1_{\rho }\times B^{n-1}_{\rho}$. Considérons maintenant une famille de formes positives 
$$
(\delta _{\varepsilon})_{\varepsilon >0}\subset {\cal D}^{n-1,n-1}(B^{n-1}_{\rho})^+
$$
convergentes faiblement vers le courant de Dirac $\delta _0\in {{\cal D}'}^{n-1,n-1}(B^{n-1}_{\rho})^+$ en $0$ lorsque $\varepsilon$ tend vers $0$, (sur le cône des courants positifs la topologie faible coïncide avec la topologie faible des courants d'ordre zéro). Si on désigne par $p_2:B^1_{\rho }\times B^{n-1}_{\rho}\longrightarrow B^{n-1}_{\rho}$ la deuxième projection on a 
$$
[B^1_{\rho }\times 0]=\lim_{\varepsilon \rightarrow 0}p_2^*\delta _{\varepsilon},
$$ 
où la limite est considéré dans la topologie faible des courants d'ordre zéro. 
\\
Pour tout forme  $\varphi \in {\cal D}^{1,1}(B^1_{\rho }\times B^{n-1}_{\rho}, J_0)^+$ positive par raport à la stucture presque complexe canonique de $\C^n$ on a
$$
0\leq \int\limits_{B^1_{\rho }\times B^{n-1}_{\rho}} \Delta_{z_1}(f\circ \sigma) \, p_2^*\delta _{\varepsilon}\wedge \varphi=
\int\limits_{B^1_{\rho }\times B^{n-1}_{\rho}}  (f\circ \sigma)\,p_2^*\delta _{\varepsilon}\wedge  \Delta_{z_1}\varphi.
$$
Si on désigne par $j:B^1_{\rho }\rightarrow B^1_{\rho }\times B^{n-1}_{\rho}, \,j(B^1_{\rho })=B^1_{\rho }\times 0$, l'immersion canonique on a
$$
\int\limits_{B^1_{\rho }} \Delta (f\circ \gamma )\,j^*\varphi =
\int\limits_{B^1_{\rho }} (f\circ \gamma )\,\Delta(j^*\varphi) =
\lim_{\varepsilon \rightarrow 0} \int\limits_{B^1_{\rho }\times B^{n-1}_{\rho}}  (f\circ \sigma)\,p_2^*\delta _{\varepsilon}\wedge  \Delta_{z_1}\varphi
\geq 0.
$$
La surjectivité de l'application $\varphi \in {\cal D}^{1,1}(B^1_{\rho }\times B^{n-1}_{\rho}, J_0)^+\mapsto j^*\varphi \in {\cal D}^{1,1}(B^1_{\rho })^+$ permet alors de conclure.\hfill $\Box$
\section{Sur la régularisation des potentiels sur les  variétés presque complexes avec contrôle asymptotique de la perte de 
\\
positivité du courant}
Pour la solution de la conjecture dans le cas général d'une distribution réelle $u$ (qui est un élément de $W^{1,1}_{loc}(X)$ d'après la remarque 2 de la section précédente) nous proposons une technique de régularisation globale des potentiels $u$ des $(1,1)$-courants positifs du type $i\partial_{_J }\bar{\partial}_{_J }u$ sur les variétés presque complexes analogue à celle utilisé avec succès par Demailly \cite{Dem-2} dans le cas complexe intégrable. La nécessité d'utiliser une technique globale dérive du fait que sur une variété presque complexe non intégrable on a pas de coordonnées naturelles qui permettent de régulariser $u$ sans perte de positivité du courant $i\partial_{_J }\bar{\partial}_{_J }u$.
\\
\\
Soit $(X,J)$ une variété presque complexe et $\omega \in {\cal E}(\Lambda ^{1,1}_{_J}T_X^*)(X)$ une métrique hermitienne sur $T_{X,J}$. Soit $\exp^{\omega} : {\cal U}\subset T_X\longrightarrow X$ le flot géodésique induit par la connexion de Chern du fibré tangent
$$
D^{\omega} _{_{J} }: {\cal E}(T_{X,J}) \longrightarrow {\cal E}(T^*_X\otimes_{_{\R}}T_{X,J})
$$
associé à la métrique $\omega$, (ici ${\cal U}\subset T_X$ désigne un voisinage ouvert de la section nulle). On désigne par
$$
\exp^{\omega} _{x,\varepsilon}:=\exp^{\omega} _x(\varepsilon \cdot):T_{X,x}\cap \varepsilon ^{-1}{\cal U}  \longrightarrow X,\,\varepsilon>0,
$$
et  on considère une fonction 
$\chi:\R\longrightarrow \R$ de classe $\ci$ telle que $\chi(t)>0$ pour $t<1,\,\chi(t)=0$ pour $t\leq 1$ et $\int_{\C^n}\,\chi(|v|^2)\,dv=1$. On introduit alors les fonctions 
$\chi_{\varepsilon}(t):=\chi(t/\varepsilon^2)/\varepsilon^{2n}$ et l'opérateur régularisant
$$
u_{\varepsilon}(x):= \int\limits_{\zeta \in T_{X,x}}
u\circ \exp^{\omega} _x\,(\zeta)\cdot \chi_{\varepsilon}(|\zeta|^2_{\omega_x})
\frac{\omega^n_{x,\zeta } }{n!}
=
\int\limits_{\zeta \in T_{X,x}}
u\circ \exp^{\omega} _{x,\varepsilon}\,(\zeta)\cdot \chi(|\zeta|^2_{\omega_x})
\frac{\omega^n_{x,\zeta } }{n!}.
$$
Étudier la conjecture revient à étudier le contrôle asymptotique de la positivité des $(1,1)$-formes $i\partial_{_J }\bar{\partial}_{_J }u_{\varepsilon}$ car si $i\partial_{_J }\bar{\partial}_{_J }u\geq 0$ alors, comme dans le théorème de Demailly qui suivra \cite{Dem-2}, la suite de fonctions $u_{\varepsilon}$ converge de façon décroissante vers la fonction $u$ lorsque $\varepsilon$ tend vers zéro. Avant de présenter le théorème de Demailly nous introduisons les définitions suivantes.
\begin{defi}
Soit $(X,J)$ une variété presque complexe et $\omega \in {\cal E}(\Lambda ^{1,1}_{_J}T_X^*)(X)$ une métriques hermitienne sur $T_{X,J}$. On appelle courbure de Griffiths inférieure du fibré hermitien $(T_{X,J} ,\omega)$ la fonction $G_{\omega }(T_{X,J}):T_X\smallsetminus 0_X\longrightarrow \R$ homogène de degré deux 
 sur les fibres de $T_X$ définie par la formule
$$
G_{\omega }(T_{X,J})_x(\xi ):=\min_{\eta\in T_{X,x}}  
\frac{{\cal C}^{\omega}_{_{{X,J} }}(\xi \otimes \eta ,\xi \otimes \eta)}{|\eta |^2_h} 
$$
où ${\cal C}^{\omega}_{_{{X,J} }}\in {\cal E}(\Herm(T^{\otimes 2}_{X,J}))(X)$ désigne la courbure de Chern de $(T_{X,J},\omega)$. Soit 
$$
\alpha:\proj_{\C}(T_{X,J})\longrightarrow \proj_{\C}(T^*_{X,J}),\;\bar{\xi}\mapsto \alpha _{\bar{\xi}}     
$$
une application fibré de classe $\ci$ telle que $\alpha _{\bar{\xi}}\oplus \C\xi =T_{X,x},\,\xi \in T_{X,x}\smallsetminus 0_x$ pour tout $x\in X$.
On appelle courbure de Griffiths $\alpha$-partielle inférieure du fibré hermitien $(T_{X,J} ,\omega)$ la fonction 
$G^{\perp_{\alpha } }_{\omega }(T_{X,J}):T_X\smallsetminus 0_X\longrightarrow \R$ homogène de degré deux sur les fibres de $T_X$ définie par la formule
$$
G^{\perp_{\alpha } }_{\omega }(T_{X,J})_x(\xi ):=\min_{\eta\in \alpha _{\bar{\xi}} }  
\frac{{\cal C}^{\omega}_{_{{X,J} }}(\xi \otimes \eta ,\xi \otimes \eta)}{|\eta |^2_h}. 
$$
\end{defi}  
On remarque que le nombre $G_{\omega }(T_{X,J})$ est la plus petite valeur propre de l'endomorphisme hermitien $i{\cal C}_{\omega}(T_{X,J})(\xi ,J\xi )$ de $(T_{X,J} ,\omega)$, où ${\cal C}_{\omega}(T_{X,J})$ désigne le tenseur de courbure de Chern de $(T_{X,J} ,\omega)$. De plus on a l'inégalité évidente $G^{\perp_{\alpha} }_{\omega }(T_{X,J})\geq G_{\omega }(T_{X,J})$. Une application $\alpha$ peut étre induite par une métrique hermitienne sur $T_X$ par exemple. Il est facile de voir (cf. \cite{Dem-2}) qu'il existe une métrique hermitienne $\omega$ sur $\proj^1_{\C}\times C$, où $C$ est une courbe holomorphe de gendre $g\geq 2$ telle que $G^{\perp_{\omega } }_{\omega }(T_{X,J})>0$. On peut énoncer le résultat de Demailly sur la forme suivante.
\begin{theo}{\bf (Demailly)} 
Soit $(X,J,\omega)$ une variété complexe hermitienne telle que 
\\
$G^{\perp_{\alpha } }_{\omega }(T_{X,J})\geq 0$ 
pour une certaine $\alpha$ et soit $\psi$ une fonction quasi-plurisousharmonique telle que $i\partial_{_J }\bar{\partial}_{_J }\psi \geq \gamma$, où $\gamma$ est une $(1,1)$-forme continue. Soit $\exph^{\omega}$ la partie holomorphe sur les fibres du fibré tangent de l'application exponentielle $\exp^{\omega}$.
L'opérateur régularisant
$$
\psi_{\varepsilon}(x):= \int\limits_{\zeta \in T_{X,x}}
\psi \circ \exph^{\omega} _x\,(\zeta)\cdot \chi_{\varepsilon}(|\zeta|^2_{\omega_x})
\frac{\omega^n_{x,\zeta } }{n!}
$$
vérifie les propriétées suivantes:
\\
{\bf A)} la suite de fonctions $(\psi_{\varepsilon} )_{\varepsilon>0}\subset {\cal E}(X,\R)$ converge ponctuellement de façon décroissante vers $\psi$ lorsque $\varepsilon>0$ tend vers zéro,  
\\
{\bf B)} on a le contrôle 
$i\partial_{_J }\bar{\partial}_{_J }\psi_{\varepsilon}  \geq \gamma-\delta _{\varepsilon}\omega$ sur la perte de positivité, où $(\delta _{\varepsilon})_{\varepsilon>0}\subset (0,+\infty)$ est une famille de 
réels qui tend vers zéro de façon décroissante lorsque $\varepsilon>0$ tend vers zéro.
\end{theo}  
Le théorème précédent reste valable même si on considère simplement l'application exponentielle $\exp^{\omega}$, mais les calculs de la preuve deviennent plus compliquées. Dans le cas presque complexe non intégrable on ne peut pas envisager de définir l'application $\exph^{\omega}$, pour des raisons délicates sur le jet d'ordre deux de la structure presque complexe.\\
\\
On estime que le lemme $\ref{Hess2}$ de la sous-section qui suivra peut être utile pour la preuve de la conjecture. On rappelle d'abord la notion de coordonnées presque complexes centrées en un point, notion qui a été introduite dans \cite{Pal}. On rappelle le corollaire suivant qui à été montré dans \cite{Pal} .
\begin{coro}\label{J3}
Pour tout point $x$ d'une variété presque complexe $(X,J)$ il existe des coordonnées $(z_1,...,z_n)$ de classe $\ci$ centrées en $x$ telles que les matrices $A(z)$ et $B(z)$ de la structure presque complexe 
$$
J(z)=\sum_{k,l}\Big(A_{k,l}(z)\,dz_l\otimes \frac{\partial}{\partial z_k} +B_{k,l}(z)\, dz_l\otimes\frac{\partial}{\partial \bar{z}_k}+\overline{B}_{k,l}(z)\, d\bar{z}_l\otimes \frac{\partial}{\partial z_k}
+\overline{A}_{k,l}(z)\, d\bar{z}_l\otimes\frac{\partial}{\partial \bar{z}_k}  \Big),
$$
relatives à ces coordonnées admettent les développements asymptotiques
\begin{eqnarray}
&B(z)=\displaystyle{\sum_r \,B^rz_r+\sum_{r,s}}\, \Big(B^{r,s}\,z_rz_s +B^{r,\bar{s}}\,z_r\bar{z}_s \Big)&\nonumber
\\\nonumber
\\
&\displaystyle{+\sum_{r,s,t}\,\Big(B^{r,s,t}\,z_rz_sz_t+B^{r,s,\bar{t} }\,z_rz_s\bar{z}_t+B^{r,\bar{s},\bar{t} }\,z_r\bar{z}_s\bar{z}_t\Big) +O(|z|^4)}&\label{B3}
\\\nonumber
\\
&A (z)=i\,I_n+\displaystyle{\frac{i}{2}\sum_{r,s} \,\overline{B}^{r} \cdot B^{s} \,z_s\bar{z}_r+
\frac{i}{4 } \sum_{r,s,t}\,\Big(\overline{B}^{t,\bar{r} }\cdot B^{s}+\overline{B}^{t,\bar{s}}\cdot B^r+2\overline{B}^t\cdot B^{r,s}\Big)\, z_rz_s\bar{z}_t  }&\nonumber
\\\nonumber
\\
&\displaystyle{+\frac{i}{4 } \sum_{r,s,t}\,\Big(\overline{B}^t\cdot B^{r,\bar{s}}+\overline{B}^s\cdot B^{r,\bar{t} }+2 \overline{B}^{s,t}\cdot B^r\Big)\, z_r\bar{z}_s\bar{z}_t +O(|z|^4)  }&
 \label{A3}
\end{eqnarray}
où $B^r,\,B^{r,s},\,B^{r,\bar{s} },\,B^{r,s,t},\,B^{r,s,\bar{t}},\,B^{r,\bar{s},\bar{t} }\in M_{n,n}(\C)$ sont des matrices telles que $B^{r,s}$ soit symétrique par rapport aux indices $r,s$, $B^{r,s,t}$ par rapport à $r,s,t$, 
$B^{r,s,\bar{t} }$ par rapport à $r,s$, $B^{r,\bar{s},\bar{t}}$ par rapport à $s,t$ et $B^r_{k,l}=0$ pour $r\leq l$, $B^{r,s}_{k,l}=0 $ pour $r,s\leq l$, $B^{r,\bar{s}}_{k,l}=0$ pour $r\leq l$, $B^{r,s,t}_{k,l}=0$ pour 
$r,s,t\leq l$, $B^{r,s,\bar{t}}_{k,l}=0$ pour $r,s\leq l$, et $B^{r,\bar{s} ,\bar{t} }_{k,l}=0$ pour $r\leq l$. De plus si on considère l'expression locale de la forme de torsion de la structure presque complexe
\begin{eqnarray*}
&\displaystyle{\tau_{_{J}}=\sum_{1\leq k<l\leq n}[\zeta _k,\zeta _l]^{0,1}_{_{J}}\otimes \zeta ^*_k\wedge\zeta ^*_l 
=\sum_{\substack{1\leq k<l\leq n
\\
1\leq r\leq n }}\overline{N}^r_{k,l}\, \zeta ^*_k\wedge\zeta ^*_l \otimes \bar{\zeta}_r }&
\end{eqnarray*}
où $\zeta _l:=(\partial /\partial z_l)^{1,0}_{_{J} }\in {\cal E}(T^{1,0}_{X,J})(U_x), \,l=1,...,n$ est le repère locale du fibré des $(1,0)$-vecteurs $T^{1,0}_{X,J}$ issue des coordonnées $(z_1,...,z_n)$ on a l'expression
$$
\overline{N}^r_{k,l}(z)=\frac{i}{2}\,B^l_{r,k}+\frac{i}{2}\sum_s\Big[2(B^{l,s}_{r,k}-B^{k,s}_{r,l})\,z_s+ B^{l,\bar{s}}_{r,k} \,\bar{z}_s\Big]+O(|z|^2)
$$
pour tout $k<l$. Le jet d'ordre $k=0,1$ de la forme de torsion de la structure presque complexe au point $x$ est nul si et seulement si les coefficients $B_{*,*}(z)$ de la structure presque complexe relatifs aux coordonnées en question s'annulent à l'ordre $k+1$.
\end{coro}
Les coordonnées précédentes sont appelées coordonnées presque complexes d'ordre 3 au point $x$.
\subsection{Expression locale normale, asymptotique à l'ordre deux du Hessien presque complexe} 
Nous avons la définition suivante.
\begin{defi}\label{DefadHess} 
Soit $(X,J)$  une variété presque complexe de dimension complexe $n$.
Pour tout distribution réelle $u\in {\cal D}'_{2n}(\R)(X)$ et tout champ de vecteurs réel $\xi \in {\cal E}(T_X)(X)$  on définit le Hessien 
presque complexe par la formule 
$$
H_{_J }u\,(\xi):=i\partial_{_J }\bar{\partial}_{_J }u\,(\xi ,J\xi )=
2\partial_{_J }\bar{\partial}_{_J }u\,(\xi^{1,0}  ,\xi^{0,1}  )=
\frac{1}{2} (\xi .\,\xi .\,u+ J\xi .\,J\xi .\,u+J[\xi ,J\xi ].\,u ).
$$
\end{defi}  
Avec les notations du corollaire $\ref{J3}$ on a le lemme  suivant.
\begin{lem}\label{Hess2} 
Soit $(X,J)$ une variété presque complexe et $u\in {\cal D}'_{2n}(\R)(X)$ une distribution réelle sur $X$. Soient $(z_1,...,z_n)$ des coordonnées locales presque complexes d'ordre $N\geq 3$ en un point $x\in X$. Alors pour tout champ de vecteurs réels 
$
\xi= \sum_k\Big(\xi _k \frac{\partial}{\partial z_k}+\bar{\xi}_k\frac{\partial}{\partial \bar{z}_k}\Big)
$
on a l'expression asymptotique à l'ordre deux du Hessien presque complexe suivante:
\begin{eqnarray*}
&\displaystyle{H_{_J }u\,(\xi)(z)=2\sum_{k,l} \frac{\partial^2u}{\partial z_k\partial\bar{z}_l}(z)\,\xi _k\bar{\xi }_l
+\sum_{k,l}\,\Re e\Big[ Q_{k,\bar{l} }(z,\xi ) \frac{\partial^2 u}{\partial z_k\partial\bar{z}_l}(z) +Q_{k,l }(z,\xi )\frac{\partial^2u}{\partial z_k\partial z_l}(z) \Big]+}
\\
\\
&\displaystyle{+\sum_{k,l,s}  \,\Re e\Big[ \Big( R^s_{k,l}(z)\,\xi _k\xi _l+R^s_{k,\bar{l}}(z)\,\xi _k\bar{\xi }_l+ R^s_{\bar{k},\bar{l}}(z)\,\bar{\xi }_k\bar{\xi }_l\Big) \frac{\partial u}{\partial z_s}(z) \Big]+O(|z|^3)(\xi ,\bar{\xi })}&
\end{eqnarray*}
où 
\begin{eqnarray*}
&\displaystyle{Q_{k,\overline{l} }(z,\xi ) :=2i\sum_t\,{\bf jet}_2B_{l,t}(z)  \,\xi _k\xi _t+ \sum_{s,t,r,h} \,\overline{B}^h_{k,s} \Big( B^r_{s,t}\,\xi _t\bar{\xi }_l }+
 B^r_{l,t}\,\xi _t\bar{\xi }_s \Big)\,z_r\bar{z}_h ,&
\\
\\
&\displaystyle{Q_{k,l }(z,\xi ) := i\sum_t\,\Big(\overline{ {\bf jet}_2B_{k,t}(z) }\,\xi _l\bar{\xi }_t+\overline{ {\bf jet}_2B_{l,t}(z) }\,\xi _k\bar{\xi }_t\Big)- }&
\\
\\
&\displaystyle{ -  \sum_{s,t,r,h}\,\Big[ \frac{1}{2}\Big(\overline{B}^h_{k,s} \,\xi _l\xi _t
+\overline{B}^h_{l,s}\,\xi _k\xi _t\Big)  B^r_{s,t}\,z_r\bar{z}_h-\overline{B}^r_{k,s} \overline{B}^h_{l,t}\,  \bar{z}_r\bar{z}_h\,\bar{\xi }_s \bar{\xi }_t  \Big]  }&  
\end{eqnarray*}
et ${\bf jet}_2B_{l,t}(z)$ est la composante $(l,t)$ du jet d'ordre deux au point zéro de la matrice $B(z)$ de la structure presque complexe par rapport aux coordonnées en question.
De plus
\begin{eqnarray*}
&\displaystyle{R^s_{k,l}(z):=\sum_{r,h}\,\Big(R^{s,r,h}_{k,l}\,z_rz_h+R^{s,r,\bar{h}}_{k,l}\,z_r\bar{z}_h+R^{s,\bar{r},\bar{h}}_{k,l}\,\bar{z}_r\bar{z}_h \Big) }& 
\\
\\
&\displaystyle{R^s_{k,\bar{l}}(z):=\sum_r\,\Big(\sum_t\,B^r_{t,k} \overline{B}^t_{s,l}\,z_r+2i\, \overline{B}_{s,l}^{r,\bar{k}}\,\bar{z}_r \Big)}&
\\
\\
&\displaystyle{
+\sum_{r,h}\,\Big(R^{s,r,h}_{k,\bar{l}}\,z_rz_h+R^{s,r,\bar{h}}_{k,\bar{l}}\,z_r\bar{z}_h+R^{s,\bar{r},\bar{h}}_{k,\bar{l}}\,\bar{z}_r\bar{z}_h \Big) }& 
\\
\\
&\displaystyle{R^s_{\bar{k},\bar{l}}(z):=\sum_{r,h}\,\Big(R^{s,r,\bar{h}}_{\bar{k},\bar{l}}\,z_r\bar{z}_h+R^{s,\bar{r},\bar{h}}_{\bar{k},\bar{l}}\,\bar{z}_r\bar{z}_h \Big) }& 
\end{eqnarray*}
$($voir l'appendice pour les expressions des coefficients $R^{*,*,*}_{*,*}$$)$. Enfin $O(|z|^3)(\xi ,\bar{\xi })$ désigne un polynôme homogène de degré deux par rapport aux variables $\xi _k,\;\bar{\xi }_k$ et à coefficients dans $m({\cal E}_0(\C))^3\cdot{\cal D}'_{2n}(\C)_0$, $($ici $m({\cal E}_0(\C))$ désigne l'ideal maximal dans l'anneau des germes des fonctions $\ci$ à valeurs complexes définis sur l'origine$)$ tel que $O(|z|^3)(\xi ,\bar{\xi })=\overline{O(|z|^3)(\xi ,\bar{\xi })} $. Si le jet d'ordre un de la forme de torsion de la structure presque complexe est nul au point $x$ alors l'expression assymptotique du Hessien presque complexe se réduit à la forme
\begin{eqnarray*}
&\displaystyle{
H_{_J }u\,(\xi)(z)=2\sum_{k,l} \frac{\partial^2u}{\partial z_k\partial\bar{z}_l}(z)\,\xi _k\bar{\xi }_l-}&
\\
\\
&\displaystyle{
-2\sum_{k,l,s,r,h}  \,\mbox{Im}  \Big[\Big(2\overline{B}^{h,\bar{r},\bar{k}}_{s,l}\,z_r\bar{z}_h+\overline{B}^{r,h,\bar{k}}_{s,l}\, \bar{z}_r\bar{z}_h \Big)\,\frac{\partial u}{\partial z_s}(z)\,\xi_k\bar{\xi }_l  \Big]+O(|z|^3)(\xi ,\bar{\xi })}&
\end{eqnarray*}
 Si de plus la structure presque complexe est intégrable alors par rapport à tout coordonnées complexes centrées en $x\in X$ on a l'expression classique
$$
H_{_J }u\,(\xi)(z)=2\sum_{k,l} \frac{\partial^2u}{\partial z_k\partial\bar{z}_l}(z)\,\xi _k\bar{\xi }_l
$$
\end{lem}  
$Preuve$.
On va écrire les expressions des opérateurs différentiels $J\xi .\,J\xi$ et $J[\xi ,J\xi ]$. Avec les notations introduites précédemment on a 
\begin{eqnarray}
J\xi=\sum_{k,l}\Big[(A_{k,l}\xi _l+\overline{B}_{k,l}\bar{\xi }_l   )\frac{\partial}{\partial z_k}+(\overline{A}_{k,l}\bar{\xi }_l+B_{k,l}\xi _l)\frac{\partial}{\partial \bar{z}_k}\Big] \label{Jx} 
\end{eqnarray}
Bien évidemment il suffit d'effectuer nos calculs pour des champs  de vecteurs réels à coefficients constants. A partir de maintenant on va donc supposer que le champ $\xi $ est à coefficients constants par rapport aux coordonnées presque complexes en considération.
On a alors $J\xi .\,J\xi=T+\overline{T} $ où
\begin{eqnarray*}
T:=\sum_{k,l,s,t}\, \Big(A_{k,l}\xi _l+\overline{B}_{k,l}\bar{\xi }_l \Big) \Big(\frac{\partial A_{s,t}}{\partial z_k}\xi _t+
\frac{\partial\overline{B}_{s,t}}{\partial z_k}\bar{\xi }_t\Big)\frac{\partial}{\partial z_s}
\\
\\
+\sum_{k,l,s,t}\,\Big(A_{k,l}\xi _l+\overline{B}_{k,l}\bar{\xi }_l  \Big )  \Big(A_{s,t}\xi _t+\overline{B}_{s,t}\bar{\xi }_t\Big )\frac{\partial^2}{\partial z_k\partial z_s}
\\
\\
+\sum_{k,l,s,t}\,\Big(A_{k,l}\xi _l+\overline{B}_{k,l}\bar{\xi }_l \Big) \Big(\frac{\partial\overline{A}_{s,t}}{\partial z_k}\bar{\xi }_t+
\frac{\partial B_{s,t}}{\partial z_k}\xi _t\Big)\frac{\partial}{\partial \bar{z}_s}
\\
\\
+\sum_{k,l,s,t}\,\Big(A_{k,l}\xi _l+\overline{B}_{k,l}\bar{\xi }_l  \Big )  \Big(\overline{A}_{s,t}\bar{\xi }_t+B_{s,t}\xi_t\Big )\frac{\partial^2}{\partial z_k\partial\bar{z}_s}
\end{eqnarray*}
En utilisant les expressions $\eqref{A3} $ et $\eqref{B3}$ on a 
\begin{eqnarray*}
&\displaystyle{\frac{\partial A}{\partial z_k}(z)=\frac{i}{2}\sum_r \,\overline{B}^{r}\cdot B^k\,\bar{z}_r +\frac{i}{2}\sum_{r,t}\,\Big(\overline{B}^{t,\bar{k} }\cdot B^r+\overline{B}^{t,\bar{r} }\cdot B^k+2\overline{B}^t \cdot  B^{k,r}  \Big)\,z_r\bar{z}_t  }&
\\
\\
&\displaystyle{+\frac{i}{4}\sum_{r,t}\,\Big(\overline{B}^{t}\cdot B^{k,\bar{r}}+\overline{B}^r\cdot B^{k,\bar{t} }+2 \overline{B}^{r,t}\cdot B^k\Big)\,\bar{z}_r\bar{z}_t  +O(|z|^3)}&
\end{eqnarray*}
\begin{eqnarray*}
&\displaystyle{\frac{\partial A}{\partial \bar{z}_k}(z)=\frac{i}{2}\sum_r \,\overline{B}^k\cdot B^r\,z_r +\frac{i}{4}\sum_{r,t}\,\Big(\overline{B}^{k,\bar{r} }\cdot B^t+ \overline{B}^{k,\bar{t} }\cdot B^r +2\overline{B}^k\cdot B^{r,t}\Big)\,z_rz_t }&
\\
\\
&\displaystyle{+\frac{i}{2}\sum_{r,t}\,\Big(\overline{B}^t\cdot B^{r,\bar{k} }+\overline{B}^k\cdot B^{r,\bar{t}}+2\overline{B}^{k,t}\cdot B^r\Big)\,z_r\bar{z}_t +O(|z|^3)  }&
\\
\\
&\displaystyle{\frac{\partial B}{\partial z_k}(z)=B^k+\sum_r}\, \Big(2B^{k,r}z_r +B^{k,\bar{r} }\bar{z}_r\Big) +\sum_{r,s}\,\Big(3B^{k,r,s}z_rz_s +2B^{k,r,\bar{s} }z_r\bar{z}_s+B^{k,\bar{r},\bar{s}}\bar{z}_r \bar{z}_s \Big)+O(|z|^3)&
\\
\\
&\displaystyle{\frac{\partial B}{\partial\bar{z}_k} (z)=\sum_r  B^{r,\bar{k} }z_r+\sum_{r,s} }\, \Big(B^{r,s,\bar{k}}z_rz_s+2B^{r,\bar{s},\bar{k}}z_r\bar{z}_s\Big)+O(|z|^3)&
\end{eqnarray*}
On explicite maintenant les quatre termes de l'opérateur $T$ à l'aide des expressions précédentes sur les dérivées premières des matrices $A$, $B$ et des relations $\eqref{A3}$, $\eqref{B3}$.
\\
I-er terme de l'opérateur $T$: 
\begin{eqnarray*}
&\displaystyle{\sum_{k,l,s,t}\, \Big(A_{k,l}\xi _l+\overline{B}_{k,l}\bar{\xi }_l \Big)\Big(\frac{\partial A_{s,t}}{\partial z_k}\xi _t+
\frac{\partial\overline{B}_{s,t}}{\partial z_k}\bar{\xi }_t\Big)\frac{\partial}{\partial z_s}=}&
\\
\\
&\displaystyle{=\sum_{k,s,t}\,i\Big(\frac{\partial A_{s,t}}{\partial z_k}\xi _k\xi _t+
\frac{\partial\overline{B}_{s,t}}{\partial z_k}\xi _k\bar{\xi }_t\Big)\frac{\partial}{\partial z_s}}&
\\
\\
&\displaystyle{+\frac{i}{2}\sum_{k,l,s,t,r,h,j}\,\Big(\overline{B}^r_{k,j}  B^h_{j,l}  B^k_{s,t}\,z_h\bar{z}_r  +\overline{B}^r_{k,t} \overline{B}^h_{s,j}  B^k_{j,l}\,\bar{z}_r\bar{z}_h\Big)\,\xi_l\bar{\xi }_t\frac{\partial}{\partial z_s} }&
\\
\\
&\displaystyle{+\sum_{k,l,s,t,r,h}\,\overline{B}^r_{k,l}   \overline{B}^{h,\bar{k} }_{s,t}\,\bar{z}_r\bar{z}_h\,\bar{\xi }_l\bar{\xi }_t\frac{\partial}{\partial z_s}
+O(|z|^3)\frac{\partial}{\partial z}  }&
\end{eqnarray*}
II-ème terme de l'opérateur $T$:
\begin{eqnarray*}
&\displaystyle{\sum_{k,l,s,t}\,\Big(A_{k,l}\xi _l+\overline{B}_{k,l}\bar{\xi }_l  \Big )  \Big(A_{s,t}\xi _t+\overline{B}_{s,t}\bar{\xi }_t\Big )\frac{\partial^2}{\partial z_k\partial z_s}=}&
\\
\\
&\displaystyle{=\sum_{k,l}\,\Big[-\xi _k\xi _l+ 2i\sum_t\,\overline{ {\bf jet}_2B_{l,t}(z)}  \,\xi _k\bar{\xi }_t-
\sum_{r,h,j,t}\, \overline{B}^r_{k,j}  B^h_{j,t}\,z_h\bar{z}_r \,\xi _l\xi _t   \Big]\frac{\partial^2}{\partial z_k\partial z_l} }&
\\
\\
&\displaystyle{+\sum_{k,l,s,t,r,h}\,\overline{B}^r_{k,s}  \overline{B}^h_{l,t}\,\bar{z}_r\bar{z}_h\,\bar{\xi }_s \bar{\xi }_t \frac{\partial^2}{\partial z_k\partial z_l}
+O(|z|^3)\frac{\partial^2}{\partial z^2} }&
\end{eqnarray*}
III-ème terme de l'opérateur $T$ :
\begin{eqnarray*}
&\displaystyle{\sum_{k,l,s,t}\,\Big(A_{k,l}\xi _l+\overline{B}_{k,l}\bar{\xi }_l \Big) \Big(\frac{\partial\overline{A}_{s,t}}{\partial z_k}\bar{\xi }_t+
\frac{\partial B_{s,t}}{\partial z_k}\xi _t\Big)\frac{\partial}{\partial \bar{z}_s}=}&
\\
\\
&=\displaystyle{\sum_{k,s,t}\,i\Big(\frac{\partial \overline{A}_{s,t}}{\partial z_k}\bar{\xi }_t\xi_k+
\frac{\partial B_{s,t}}{\partial z_k}\xi _k\xi _t\Big)\frac{\partial}{\partial \bar{z}_s} }&
\end{eqnarray*}
\begin{eqnarray*}
&\displaystyle{+\sum_{k,l,s,t,r}\,\Big[ \overline{B}^r_{k,l}  B^k_{s,t}\,\bar{z}_r\,\bar{\xi }_l\xi _t
+\frac{i}{2}\sum_{j,h}\,\Big(\overline{B}^r_{k,j}  B^h_{j,l}  B^k_{s,t}z_h\bar{z}_r\, \xi _l\xi _t- \overline{B}^r_{k,l}  B^k_{s,j}   \overline{B}^h_{j,t}\bar{z}_h \bar{z}_r\,\bar{\xi }_l\bar{\xi }_t\Big)\Big]\frac{\partial}{\partial \bar{z}_s} }&
\\
\\
&\displaystyle{+\sum_{k,l,s,t,r,h}\,\Big[ \Big(2\overline{B}^r_{k,l}  B^{k,h}_{s,t}+\overline{B}^{r,\bar{h} }_{k,l}   B^k _{s,t}\Big)\bar{z}_rz_h +\Big(\overline{B}^r_{k,l}  B^{k,\bar{h} }_{s,t} + \overline{B}^{r,h}_{k,l}  B^k _{s,t}\Big) \bar{z}_r\bar{z}_h  \Big]\,\bar{\xi }_l\xi _t \frac{\partial}{\partial \bar{z}_s}+O(|z|^3)\frac{\partial}{\partial \bar{z}}}&
\end{eqnarray*}
IV-ème terme de l'opérateur $T$:
\begin{eqnarray*}
&\displaystyle{\sum_{k,l,s,t}\,\Big(A_{k,l}\xi _l+\overline{B}_{k,l}\bar{\xi }_l  \Big )  \Big(\overline{A}_{s,t}\bar{\xi }_t+B_{s,t}\xi_t\Big )\frac{\partial^2}{\partial z_k\partial\bar{z}_s}=}&
\\
\\
&\displaystyle{=\sum_{k,l}\,\xi _k\bar{\xi }_l\frac{\partial^2}{\partial z_k\partial\bar{z}_l}  
+\frac{1}{2}  \sum_{k,l,t,r,h,j}\,\Big(\overline{B}^r_{k,j}  B^h_{j,t}\,z_h\bar{z}_r\,\xi _t\bar{\xi }_l\frac{\partial^2}{\partial z_k\partial\bar{z}_l} 
+B^r_{k,j}  \overline{B}^h_{j,t}\,\bar{z}_hz_r\,\bar{\xi } _t\xi _l\frac{\partial^2}{\partial \bar{z}_k\partial z_l}\Big) }&
\\
\\
&\displaystyle{+i\sum_{k,l,t}\,\Big({\bf jet}_2B_{l,t}(z)\,\xi _k\xi _t\frac{\partial^2}{\partial z_k\partial\bar{z}_l}
-\overline{ {\bf jet}_2B_{l,t}(z) }\,\bar{\xi }_k\bar{\xi }_t\frac{\partial^2}{\partial\bar{z}_k\partial z_l } \Big) } & 
\\
\\
&\displaystyle{+\sum_{k,l,s,t,r,h}\,\overline{B}^r_{k,s}  B^h_{l,t}\,\bar{z} _rz_h\, \bar{\xi }_s\xi _t\frac{\partial^2}{\partial z_k\partial\bar{z}_l}    +O(|z|^3)\frac{\partial}{\partial z\partial \bar{z}}  }&
\end{eqnarray*}
 Venons-en maintenant au calcul du champ de vecteurs $J[\xi ,J\xi]$. On remarque d'abord que:
\begin{eqnarray*}
[\xi ,J\xi ]=\sum_{k,s,t}\,\Big(\frac{\partial A_{s,t}}{\partial z_k}\xi _k\xi _t+
\frac{\partial\overline{B}_{s,t}}{\partial z_k}\xi _k\bar{\xi }_t+\frac{\partial A_{s,t}}{\partial \bar{z}_k}\xi _t\bar{\xi }_k+
\frac{\partial\overline{B}_{s,t}}{\partial\bar{z}_k}\bar{\xi }_k\bar{\xi }_t\Big)\frac{\partial}{\partial z_s}
\\
\\
+\sum_{k,s,t}\, \Big(\frac{\partial\overline{A}_{s,t}}{\partial \bar{z}_k}\bar{\xi }_k\bar{\xi }_t+
\frac{\partial B_{s,t}}{\partial \bar{z}_k}\bar{\xi }_k\xi _t+\frac{\partial \overline{A}_{s,t}}{\partial z_k}\bar{\xi }_t\xi_k+
\frac{\partial B_{s,t}}{\partial z_k}\xi _k\xi _t\Big)\frac{\partial}{\partial \bar{z}_s}
\end{eqnarray*}
En appliquant la relation $\eqref{Jx}$ avec $[\xi ,J\xi ]$ à la place de $\xi $ on a $J[\xi ,J\xi ]=\eta+\bar{\eta}$ avec:
\begin{eqnarray*}
\eta:=\sum_{k,s,t,r}\,\Big[A_{r,s}  \Big(\frac{\partial A_{s,t}}{\partial z_k}\xi _k\xi _t+
\frac{\partial\overline{B}_{s,t}}{\partial z_k}\xi _k\bar{\xi }_t+\frac{\partial A_{s,t}}{\partial \bar{z}_k}\xi _t\bar{\xi }_k+
\frac{\partial\overline{B}_{s,t}}{\partial\bar{z}_k}\bar{\xi }_k\bar{\xi }_t\Big)
\\
\\
+\overline{B}_{r,s} \Big(\frac{\partial\overline{A}_{s,t}}{\partial \bar{z}_k}\bar{\xi }_k\bar{\xi }_t+
\frac{\partial B_{s,t}}{\partial \bar{z}_k}\bar{\xi }_k\xi _t+\frac{\partial \overline{A}_{s,t}}{\partial z_k}\bar{\xi }_t\xi_k+
\frac{\partial B_{s,t}}{\partial z_k}\xi _k\xi _t\Big)  \Big]\frac{\partial}{\partial z_r}
\end{eqnarray*}
En tenant compte des relations $\eqref{A3}$, $\eqref{B3}$ et celles sur les dérivées premières des matrices $A$, $B$ dans l'expression du champ $\eta$ on obtient:
\begin{eqnarray*}
&\displaystyle{\eta=i\sum_{k,s,t}\, \Big(\frac{\partial A_{s,t}}{\partial z_k}\xi _k\xi _t+
\frac{\partial\overline{B}_{s,t}}{\partial z_k}\xi _k\bar{\xi }_t+\frac{\partial A_{s,t}}{\partial \bar{z}_k}\xi _t\bar{\xi }_k+
\frac{\partial\overline{B}_{s,t}}{\partial\bar{z}_k}\bar{\xi }_k\bar{\xi }_t\Big)\frac{\partial}{\partial z_s} }&
\\
\\
&\displaystyle{+\sum_{k,l,s,t,r,h}\,\Big(\overline{B}^l_{r,s}  B^{h,\bar{t}} _{s,k}\,z_h \bar{z}_l -
\frac{i}{2}\sum_j\,\overline{B}^l_{r,j}  B^k_{j,s}  \overline{B}^h_{s,t}\, \bar{z}_l\bar{z}_h\Big)\,\xi _k\bar{\xi }_t \frac{\partial}{\partial z_r} }&
\\
\\ 
&\displaystyle{+\sum_{k,l,s,t,r}\,\Big[\overline{B}^l_{r,s}  B^k_{s,t}\,\bar{z}_l+\sum_h\,\Big(2\overline{B}^l_{r,s}  B^{k,h}_{s,t}\,z_h\bar{z}_l+\overline{B}^l_{r,s}  B^{k,\bar{h}}_{s,t}\,\bar{z}_l \bar{z}_h   \Big)\Big]\,\xi _k\xi _t \frac{\partial}{\partial z_r}+O(|z|^3)\frac{\partial}{\partial z}  }&
\end{eqnarray*}
En regroupant et en simplifiant les termes obtenus jusqu'ici (et en symetrisant de façon adequate les coefficients) on obtient l'expression asymptotique voulue pour la distribution $H_{_J }u\,(\xi)$. Si maintenant on suppose que la structure presque complexe est intégrable il existe d'après le théorème de Newlander-Nirenberg des coordonnées locales complexes. La structure presque complexe s'écrit alors par rapport à ces coordonnées sous la forme 
$$
J(z)=J_0=i\sum_k\Big(dz_k\otimes \frac{\partial}{\partial z_k}-d\bar{z}_k\otimes\frac{\partial}{\partial \bar{z}_k}\Big).
$$
Un calcul immédiat montre alors que pour tout champs de vecteurs réels $\xi $ à coefficients constants par rapport aux coordonnées complexes on a: 
$$
\frac{1}{2}(\xi .\xi +J\xi .J\xi )=2\sum_{k,l}\,\xi _k\bar{\xi }_l  \frac{\partial^2}{\partial z_k\partial\bar{z}_l}, 
$$ 
(bien évidemment $[\xi ,J\xi ]=0$), ce qui permet de conclure dans le cas  d'une structure presque complexe intégrable. \hfill $\Box$

\section{Appendice} 
\subsection{Expressions des coefficients $R^{*,*,*}_{*,*}$ du Hessien presque complexe} 
\begin{eqnarray*}
&\displaystyle{R_{k,l}^{s,r,h}:=\frac{i}{8}\sum_{t,j}\,\overline{B}^t_{s,j} \Big(B^r_{t,k}  B^h_{j,l}+B^h_{t,k}  B^r_{j,l}  + B^r_{t,l}  B_{j,k}^h+B^h_{t,l}  B^r_{j,k}  \Big)   }&
\\
\\
&\displaystyle{R_{k,l}^{s,r,\bar{h}}:=-\frac{1}{2}\sum_t\,\Big[\overline{B}^{h,\bar{k} }_{s,t}  B^r_{t,l}+\overline{B}^{h,\bar{l} }_{s,t}  B^r_{t,k} + \overline{B}^{h,\bar{r}}_{s,t}  \Big(B^k_{t,l}+B^l_{t,k}  \Big)\Big]  }&
\\
\\
&\displaystyle{R_{k,l}^{s,\bar{r},\bar{h}}:=-\frac{1}{2}\sum_t\,\overline{B}^{r,h}_{s,t}  \Big(B^k_{t,l}+B^l_{t,k}\Big)   }&
\\
\\
&\displaystyle{R_{k,\bar{l}}^{s,r,h}:=\frac{1}{2}\sum_t\,\Big(B^r_{t,k}  \overline{B}^{t,\bar{h}}_{s,l}+B^h_{t,k}   \overline{B}^{t,\bar{r} }_{s,l}+2B_{t,k}^{r,h}  \overline{B}^t_{s,l}     \Big)}&
\\
\\
&\displaystyle{R_{k,\bar{l}}^{s,r,\bar{h}}:=4i\,\overline{B}_{s,l}^{h,\bar{r}, \bar{k}}+\sum_t\,\Big(2B_{t,k}^r \overline{B}_{s,l}^{t,h}+B_{t,k}^{r,\bar{h}}\,  \overline{B}^t_{s,l}+  
 B_{t,k}^{r,\bar{l}} \, \overline{B}^h_{s,t}  + \frac{i}{2}\sum_j\,\overline{B}^h_{t,j}  B^r_{j,k}  B^t_{s,l}  \Big)  }&
\\
\\
&\displaystyle{R_{k,\bar{l}}^{s,\bar{r},\bar{h}}:=2i\,\overline{B}_{s,l}^{r,h,\bar{k} }+ \frac{i}{4}\sum_{t,j}\,\Big(\overline{B}^r_{t,l}  \overline{B}^h_{s,j}+\overline{B}^h_{t,l}  \overline{B}^r_{s,j}\Big) \Big( B^t_{j,k}-B^k_{j,t} \Big)  }&
\end{eqnarray*}
\begin{eqnarray*}
&\displaystyle{R_{\bar{k},\bar{l}}^{s,r,\bar{h}}:=-\frac{i}{4}\sum_{t,j}\,B^r_{t,j} \Big(\overline{B}^h_{j,k} \overline{B}^t_{s,l}+ \overline{B}^h_{j,l} \overline{B}^t_{s,k}\Big)}& 
\\
\\
&\displaystyle{R_{\bar{k},\bar{l}}^{s,\bar{r},\bar{h}}:=\frac{1}{4}\sum_t\,\Big(\overline{B}^r_{t,k} \overline{B}_{s,l}^{h,\bar{t} }+ \overline{B}^h_{t,k} \overline{B}_{s, l}^{r,\bar{t} }+
\overline{B}^r_{t,l} \overline{B}_{s,k}^{h,\bar{t} }+ \overline{B}^h_{t,l} \overline{B}_{s,k}^{r,\bar{t} }   \Big)  } .&  
\end{eqnarray*}

\vspace{1cm}
\noindent
Nefton Pali
\\
Institut Fourier, UMR 5582, Université Joseph Fourier
\\
BP 74, 38402 St-Martin-d'Hères cedex, France
\\
E-mail: \textit{nefton.pali@ujf-grenoble.fr}
\end{document}